\documentclass[preprint,12pt]{elsarticle}
\usepackage{graphics}
\usepackage{graphicx}
\usepackage{subfigure}

\usepackage{latexsym}
\usepackage{amsmath}
\usepackage{graphicx}
\usepackage{times}
\usepackage{graphicx,multicol,enumerate,subfigure,amsmath}

\usepackage{fullpage}
\usepackage{setspace}
\usepackage[english]{babel}
\usepackage[version=3]{mhchem}
\usepackage{amsmath, amssymb, amsthm}
\usepackage{verbatim, latexsym}
\usepackage{colortbl}
\usepackage{multirow}

\newcommand{\DIV}{\mbox{div}}
\newcommand{\norm}[1]{\left\|#1\right\|}

\newtheorem{theorem}{Theorem}
\newtheorem{remark}[theorem]{Remark}

\usepackage{soul}
\setstcolor{red}

\biboptions{sort&compress}
\journal{Communications in Computational Physics}

\onehalfspace

\begin{document}

\begin{frontmatter}

\title{Generalized Multiscale Finite Element Methods. Nonlinear Elliptic Equations}
\author{{Yalchin Efendiev}$^{1,2}$}

\author{{{Juan Galvis$^{3}$}, Guanglian Li$^{2}$}, and Michael Presho$^{2*}$}

\address{$^{1}$ Center for Numerical Porous Media (NumPor) \\
King Abdullah University of Science and Technology (KAUST) \\
Thuwal 23955-6900, Kingdom of Saudi Arabia.}

\address{$^{2}$ Department of Mathematics \& Institute for Scientific Computation (ISC) \\
Texas A\&M University \\
College Station, Texas, USA}

\address{$^{3}$ Departamento de Matem\'{a}ticas\\ Universidad Nacional de Colombia\\ 
Bogot\'a D.C., Colombia \\
}

\cortext[cor1]{Email address: mpresho@math.tamu.edu}

\begin{abstract}

In this paper we use the Generalized Multiscale Finite Element Method (GMsFEM) framework, introduced in \cite{egh12}, in order to solve nonlinear elliptic equations with high-contrast coefficients. The proposed solution method involves linearizing the equation so that coarse-grid quantities of previous solution iterates can be regarded as auxiliary parameters within the problem formulation. With this convention, we systematically construct respective coarse solution spaces that lend themselves to either continuous Galerkin (CG) or discontinuous Galerkin (DG) global formulations. Here, we use Symmetric Interior Penalty Discontinuous Galerkin approach. Both methods yield a predictable error decline that depends on the respective coarse space dimension, and we illustrate the effectiveness of the CG and DG formulations by offering a variety of numerical examples.

\end{abstract}


\end{frontmatter}

\section{Introduction}

Nonlinear partial differential equations represent a class of
problems that have applications in many scientific communities
\cite{egkl12,richards31}. Forchheimer flow and nonlinear
elasticity are two particular examples of physical processes
that are modeled by nonlinear equations \cite{DPG93,ga94}. In
addition to difficulties associated with the nonlinearity, these
types of problems often involve coefficients that exhibit
high-contrast, heterogeneous behavior. For example, when modeling
subsurface flow, the underlying permeability field is often represented by a
high-contrast coefficient in the pressure equation. One approach for
solving a high-contrast, nonlinear equation is to linearize the
problem and use an iterative method for obtaining the solution. For
example, a Picard iteration yields an iterative process where a
previous solution iterate is directly used in order to update the
solution at the current iteration. In this case, a final solution is
obtained when a suitable tolerance between the current and previous
iteration is reached. While relatively easy to implement, iterative
techniques typically require a repeated number of solves in order to
obtain a convergent solution. In the case of a nonlinear elliptic
equation, each iteration requires the numerical solution of a large
system of equations that depends on the previous iterate. Thus,
computing solutions on a fully resolved mesh quickly becomes a
prohibitively expensive task. As such, techniques that allow for a
more efficient computational procedure with a suitable level of accuracy
 are desirable.

 The past few decades have seen
the development of various multiscale solution techniques for
capturing small scale effects on a coarse grid
\cite{aarnes04,apwy07,eh09,hw97,hughes98,jennylt03}. The
multiscale finite element methods (MsFEM's) that we consider in this
paper hinge on the construction of coarse spaces that are spanned by
a set of independently computed multiscale basis functions. The
multiscale basis functions are then coupled via a respective global
formulation in order to compute the solution. In particular,
solutions may be computed on a coarse grid while maintaining
the fine-scale effects that are embedded into the basis functions.
While standard multiscale methods have proven effective for a
variety of applications (see, e.g.,
\cite{eghe05,eh09,ehg04,jennylt03}), in this paper we
consider a more recent framework in which the coarse spaces may be
systematically enriched to converge to the fine grid
solution \cite{bl11,egt11,egw10,nguyen08}. More specifically,
additional basis functions are chosen based on localized eigenvalue
problems that capture the underlying behavior of the system.
In this case, we may carefully choose the number of basis functions
(and dimension of the coarse space) such that we achieve a desired level of
accuracy. In this paper we additionally show that the systematic
enrichment of coarse spaces
is flexible with respect to the global formulation that is chosen to
couple the resulting basis functions.

To treat the nonlinear elliptic equation considered in this paper we
make use of the Generalized Multiscale Finite Element Method
(GMsFEM) which was introduced in \cite{egh12}.
In order to do so, we apply a Picard
iteration and treat an upscaled quantity of a previous solution
iterate as a parameter in the problem. With this convention we
follow an offline-online procedure in which the coarse space
construction is split into two distinct stages; offline and online
(see~\cite{barrault,boyoval08,nguyen08, Patera}). The main goal of
this approach is to allow for the efficient construction of an
online space (and an online solution) for each fixed parameter value
and iteration. In the process, we precompute a larger-dimensional,
\emph{parameter-independent} offline space that accounts for an
appropriate range of parameter values that may be used in the online
stage. As construction of the offline space will constitute a
one-time preprocessing step, only the online space will require
additional work within the solution procedure. In the offline stage
we first choose a fixed set of parameter values and generate an
associated set of ``snapshot" functions by solving localized
problems on specified coarse subdomains. The functions obtained
through this step constitute a snapshot space which will be used in
the offline space construction. To construct the
offline space we solve localized eigenvalue problems that use
averaged quantities of the parameter(s) of interest within the space
of snapshots. We then keep a certain number of eigenfunctions (based
on some criterion) to form the offline space. At the online stage we
solve similar localized problems using a fixed parameter value
within the offline space, and keep a certain number of
eigenfunctions for the online space construction.

In this paper we consider the continuous Galerkin (CG) and
discontinuous Galerkin (DG) formulations for the global coupling of
the online basis functions. We show that each method offers a
suitable solution technique, however, at this point we highlight some distinguishing characteristics of the
repsective methods as motivation for considering both formulations.
For the nonlinear elliptic equation considered in this paper, the CG coupling yields a bilinear form
that closely resembles the standard finite element method (FEM). In
particular, the integrations that define the CG formulation are
taken over the whole domain, and result in a reduced-order system of equations that
is similar in nature to the fine-scale system. As such, the
ease of implementation, classical FEM analogues, and well understood
structure make CG a tractable method for coupling the coarse
basis functions in order to solve the global problem \cite{hw97}. While the discontinuous Galerkin
formulation is arguably more delicate than its CG counterpart, DG
offers an attractive feature such as it does not require partition
of unity functions to couple basis functions.
Both
methods are shown to be suitable coupling mechanisms within the
GMsFEM framework that is described in this paper. In particular, an increase
in the size of the online coarse space yields a predictable error
decline, and the error is shown to behave according to previous
error estimates that depend on the eigenvalue behavior. The flexibility of the coarse space enrichment,
along with the choice of using CG or DG as the global coupling
mechanism, makes GMsFEM a robust and suitable technique for
solving the model equation that we consider in this paper. A
variety of numerical examples are presented to validate the
performance of the proposed method.

We note that some numerical results for GMsFEM in the context of continuous
Galerkin methods for nonlinear equations
are presented in \cite{egh12}. These numerical results are mostly presented
to demonstrate the main concepts of GMsFEM and we do not have careful studies
for nonlinear problems in \cite{egh12}. Moreover, the numerical results
presented in \cite{egh12} use reduced basis approach to identify dominant
eigenmodes which is different from the
local mode decomposition approach presented here.
Moreover, the current paper also studies DG approach for nonlinear equations.

 The organization of the paper is
as follows. In Sect.~\ref{prelim} we introduce the model problem,
the iterative procedure, and notation to be used throughout the
paper. In Sect.~\ref{cgdgmsfem} we carefully describe the coarse
space enrichment procedure, and introduce the continuous and
discontinuous Galerkin global coupling formulations. In particular,
Subsect.~\ref{locbasis} is devoted to the offline-online coarse
space construction, and in Subsect.~\ref{globcoupling} we describe
the CG and DG global coupling procedures. A variety of numerical
examples are presented in Sect.~\ref{numerical} to validate the
performance of the proposed approaches, and in
Sect.~\ref{conclusion} we offer some concluding remarks.

\section{Preliminaries}
\label{prelim}
In this paper we consider non-linear, elliptic equations of the form
\begin{equation} \label{eq:original}
-\mbox{div} \big( \kappa(x;u) \, \nabla u  \big)=f \, \, \text{in} \, D,
\end{equation}
where $u=0$ on $\partial D$.
We assume that $u$ is bounded above and below, i.e., $u_0\leq u(x) \leq u_N$,
where $u_0$ and $u_N$ are pre-defined constants.
We will also assume that the interval $[u_0, u_N]$ is divided into $N$ equal regions
whose endpoints are given by $u_0<u_1<...<u_{N-1}<u_N$.

\begin{figure}[htb]
  \centering
  \includegraphics[width=0.65 \textwidth]{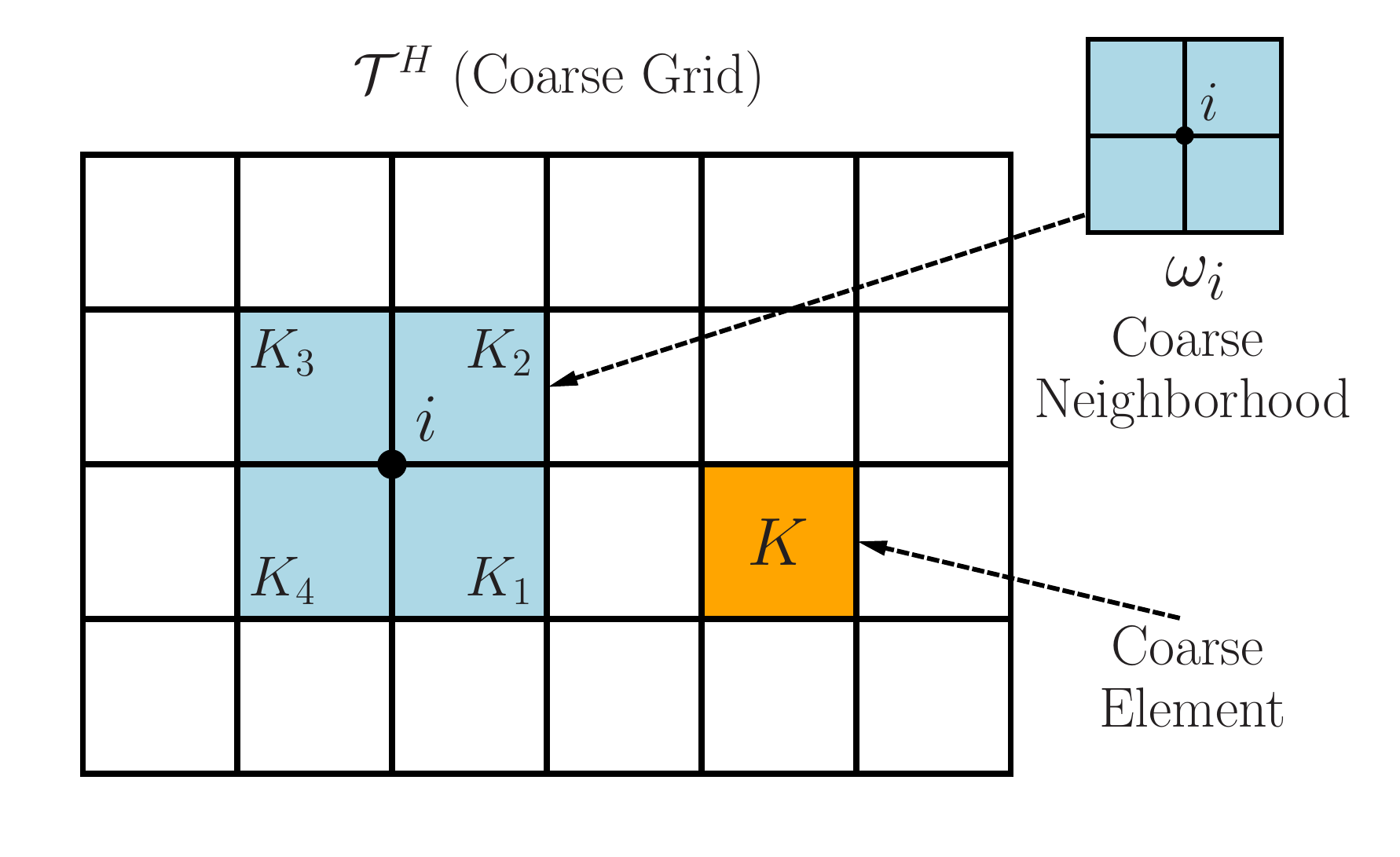}
  \caption{Illustration of a coarse neighborhood and coarse element}
  \label{schematic}
\end{figure}

In order to solve Eq.~(\ref{eq:original}) we will consider a Picard iteration
\begin{equation}
\label{eq:original1}
-\mbox{div} \big( \kappa(x;u^n(x)) \, \nabla u^{n+1}(x) \big)=f\ \text{in}\ D,
\end{equation}
where superscripts involving $n$ denote respective iteration levels. To discretize (\ref{eq:original1}), we next introduce the notion of fine and coarse grids. We let $\mathcal{T}^H$ be a usual conforming partition of the computational domain $D$ into finite elements (triangles, quadrilaterals, tetrahedrals, etc.). We refer to this partition as the coarse grid and assume that each coarse subregion is partitioned into a connected union of fine grid blocks. The fine grid partition will be denoted by $\mathcal{T}^h$. We use $\{x_i\}_{i=1}^{N_v}$ (where $N_v$ the number of coarse nodes) to denote the vertices of
the coarse mesh $\mathcal{T}^H$, and define the neighborhood of the node $x_i$ by
\begin{equation} \label{neighborhood}
\omega_i=\bigcup\{ K_j\in\mathcal{T}^H; ~~~ x_i\in \overline{K}_j\}.
\end{equation}
See Fig.~\ref{schematic} for an illustration of neighborhoods and elements subordinated to the coarse discretization. We emphasize the use of $\omega_i$ to denote a coarse neighborhood, and $K$ to denote a coarse element throughout the paper.

Next, we briefly outline the global coupling and the role of coarse basis
functions for the respective formulations that we consider. For the discontinuous Galerkin (DG) formulation, we will use a coarse element $K$ as the support for basis
functions, and for the continuous Galerkin (CG) formulation, we will use $\omega_i$ as the support of basis functions. To further motivate the coarse basis construction, we offer a brief outline of the global coupling associated with the CG formulation below. For the purposes of this description, we formally denote the CG basis functions by $\psi_k^{\omega_i}$. In particular, we note that the proposed approach will employ the use of multiple basis functions per coarse neighborhood.
In turn, the CG solution at $n$-th iteration will be sought as $u^{\text{CG}}_{\text{ms}}(x;\mu)=\sum_{i,k} c_{k}^i \psi_{k}^{\omega_i}(x; \mu)$, where $\psi_{k}^{\omega_i}(x; \mu)$ are the basis functions for $n$-th iteration, and $\mu$ is used to denote dependence on the previous solution.
We note that a main consideration of our method is to allow for rapid calculations of basis functions at each iteration.

Once the basis functions are identified, the CG global coupling is given through the variational form
\begin{equation}
\label{eq:globalG} a(u^{\text{CG}}_{\text{ms}},v;\mu)=(f,v), \quad \text{for all} \, \, v\in
V_{\text{on}}^{\text{CG}},
\end{equation}
where  $V_{\text{on}}^{\text{CG}}$ is used to denote the space formed by those basis functions.

We also note that an appropriate set of basis functions defined on
each coarse element $K$ may be respectively coupled via a
discontinuous Galerkin formulation (see e.g., \cite{MR1885715,
MR2431403,MR2002258}).

\section{CG and DG GMsFEM for nonlinear problems}
\label{cgdgmsfem}

\subsection{Local basis functions}
\label{locbasis}
To motivate the local basis construction, we first introduce
an approximation to the solution of Eq.~\eqref{eq:original1} given by
\begin{equation}
\label{eq:original2}
-\mbox{div} \big( \kappa(x;\overline{u}^n(x))\nabla u^{n+1}(x) \big) = f  \, \,  \text{in} \, D,
\end{equation}
where $\overline{u}$ denotes the average of $u$ in each coarse region
(either $K$ or $\omega_i$, depending on the desired formulation). Because the variation in
$\overline{u}^n$ is not known a priori, we use $\mu$ to represent the dependence of the solution on $\overline{u}^n$. As part of the iterative solution process, multiscale basis functions will be computed for a selected number of the parameter values at the offline stage, and we will compute
multiscale basis functions for each new value of $\overline{u}^n$ at the online stage.
In this section we will describe these details, and note that we maintain the convention of denoting
 $\overline{u}$ by the parameter $\mu$. We omit the iterative index $n$ (and $n+1$) for additional notational brevity, although note that the iterative process should be clearly implied.

With the notational conventions in place we now describe the offline-online computational procedure, and elaborate on some applicable choices for the associated bilinear forms to be used in the coarse space construction. Below we offer a general outline for the procedure.

\begin{itemize}
\item[1.]  Offline computations:
\begin{itemize}
\item 1.0. Coarse grid generation.
\item 1.1. Construction of snapshot space that will be used to compute an offline space.
\item 1.2. Construction of a small dimensional offline space by performing dimension reduction in the space of global snapshots.
\end{itemize}
\item[2.] Online computations:
\begin{itemize}
\item 2.1. For each input parameter, compute multiscale basis functions.
\item 2.2. Solution of a coarse-grid problem for any force term and boundary condition.
\item 2.3. Iterative solvers, if needed.
\end{itemize}
\end{itemize}

In the offline computation, we first construct a snapshot space $V_{\text{snap}}^{\tau}$, where $\tau$ denotes either a coarse neighborhood $\omega_i$ in the continuous Galerkin case, or a coarse element $K$ in the discontinuous Galerkin case (refer back to Fig.~\ref{schematic}). Construction of the snapshot space involves solving the local problems for various choices of input parameters, and we describe the details below.

In order to construct the space of snapshots we propose to solve the following eigenvalue problem on a coarse domain $\tau$:
\begin{equation} \label{snaplinalg}
A(\mu_j) \psi_{l,j}^{\tau, \text{snap}} = \lambda_{l,j}^{\tau, \text{snap}} S(\mu_j) \psi_{l,j}^{\tau, \text{snap}}
\quad \text{in} \, \, \, \tau,
\end{equation}
where $\mu_j$ ($j=1,\dots,J$) is a specified set of fixed parameter values, and we again emphasize that $\tau$ denotes a different coarse subdomain (either a coarse neighborhood $\omega_i$ or coarse element $K$) depending on whether we consider the CG or DG problem formulation. We are careful to note that zero Neumann boundary conditions are generally used to solve eigenvalue problem, except in the DG case when Dirichlet conditions are used on element boundaries that coincide with the global domain. The matrices in Eq.~\eqref{snaplinalg} are defined as
\begin{equation}\label{eqn:eigenmatrix}
A(\mu_j) = [a(\mu_j)_{mn}] = \int_{\tau} \kappa(x; \mu_j) \nabla \phi_n \cdot \nabla \phi_m \quad \text{and}
\quad S(\mu_j) = [s(\mu_j)_{mn}] = \int_{\tau} \widetilde{\kappa}(x; \mu_j)  \phi_n \phi_m,
\end{equation}
where $\phi_n$ denotes the standard bilinear, fine-scale basis
functions and $\widetilde{\kappa}$ will be carefully introduced in the next section. We note that Eq.~\eqref{snaplinalg} is the discretized
form of the continuous equation
\begin{equation*}
-\text{div}(\kappa(x, \mu_j) \nabla \psi_{l,j}^{\tau, \text{snap}} ) = \lambda_{l,j}^{\tau, \text{snap}} \psi_{l,j}^{\tau, \text{snap}}
\quad \text{in} \, \, \, \tau.
\end{equation*}

For brevity of notation we now omit the superscript $\tau$ for eigenvalue problems, yet it is assumed throughout this section that the offline and online space computations are localized to respective coarse subdomains. After solving Eq.~\eqref{snaplinalg}, we keep the first $L_i$ eigenfunctions corresponding to the dominant eigenvalues (asymptotically vanishing in this case) to form the space
$$
V_{\text{snap}} = \text{span}\{ \psi_{l,j}^{ \text{snap}}:~~~1\leq j \leq J ~~ \text{and} ~~ 1\leq l \leq L_i \},
$$
for each coarse neighborhood $\omega_i$ (or coarse element $K$).

 We reorder the snapshot functions using a single index to create the matrix
$$
R_{\text{snap}} = \left[ \psi_{1}^{\text{snap}}, \ldots, \psi_{M_{\text{snap}}}^{\text{snap}} \right],
$$
where $M_{\text{snap}}$ denotes the total number of functions to keep in the snapshot matrix construction.

In order to construct the offline space $V_{\text{off}}^\tau$, we perform a dimension reduction of the space of snapshots using an auxiliary spectral decomposition. The main objective is to use the offline space to efficiently (and accurately) construct a set of multiscale basis functions for each $\mu$ value in the online stage. More precisely, we seek a subspace of the snapshot space such that it can approximate any element of the snapshot space in the appropriate sense defined via auxiliary bilinear forms. At the offline stage the bilinear forms are chosen to be \emph{parameter-independent}, such that there is no need to reconstruct the offline space for each $\mu$ value. The analysis in \cite{egw10} motivates the following eigenvalue problem in the space of snapshots:
\begin{equation} \label{offeig}
A^{\text{off}} \Psi_k^{\text{off}} = \lambda_k^{\text{off}} S^{\text{off}} \Psi_k^{\text{off}},
\end{equation}
where
\begin{equation*}
 \displaystyle A^{\text{off}} = [a_{mn}^{\text{off}}] = \int_{\tau} \kappa(x, \overline{\mu}) \nabla \psi_m^{\text{snap}} \cdot \nabla \psi_n^{\text{snap}} = R_{\text{snap}}^T \overline{A} R_{\text{snap}}
 \end{equation*}
 \begin{center}
 and
 \end{center}
 \begin{equation*}
 \displaystyle S^{\text{off}} = [s_{mn}^{\text{off}}] = \int_\tau  \widetilde{\kappa}(x, \overline{\mu})\psi_m^{\text{snap}} \psi_n^{\text{snap}} = R_{\text{snap}}^T \overline{S} R_{\text{snap}},
\end{equation*}
where $\kappa(x, \overline{\mu}) $, and $\widetilde{\kappa}(x, \overline{\mu})$ are domain-based averaged coefficients with $\overline{\mu}$ chosen
as the average of pre-selected $\mu_i$'s.
 We note that $\overline{A}$ and $\overline{S}$ denote analogous fine scale matrices as defined in Eq.~\eqref{snaplinalg}, except that averaged coefficients are used in the construction. To generate the offline space we then choose the smallest $M_{\text{off}}$ eigenvalues from Eq.~\eqref{offeig} and form the corresponding eigenvectors in the space of snapshots by setting
$\psi_k^{\text{off}} = \sum_j \Psi_{kj}^{\text{off}} \psi_j^{\text{snap}}$ (for $k=1,\ldots, M_{\text{off}}$), where $\Psi_{kj}^{\text{off}}$ are the coordinates of the vector $\Psi_{k}^{\text{off}}$. We then create the offline matrix $$
R_{\text{off}} = \left[ \psi_{1}^{\text{off}}, \ldots, \psi_{M_{\text{off}}}^{\text{off}} \right]
$$
to be used in the online space construction.

For a given input parameter, we next construct the associated online coarse space
$V^{\tau}_{\text{on}}(\mu)$ \emph{for each} $\mu$ value on each coarse subdomain.
In principle, we want this to be a small dimensional subspace of the offline space for computational efficiency.
The online coarse space will be used within the finite element
framework to solve the original global problem, where a continuous or discontinuous Galerkin coupling of the multiscale basis functions is used to compute the global solution. In particular, we seek a subspace of the offline space such that it can approximate any element of the offline space in
an appropriate sense. We note that at the online stage, the bilinear forms are chosen to be \emph{parameter-dependent}. Similar analysis (see \cite{egw10}) motivates the following eigenvalue problem in the offline space:
\begin{equation} \label{oneig}
A^{\text{on}}(\mu) \Psi_k^{\text{on}} = \lambda_k^{\text{on}} S^{\text{on}}(\mu) \Psi_k^{\text{on}},
\end{equation}
where
\begin{equation*}
 \displaystyle A^{\text{on}}(\mu) = [a^{\text{on}}(\mu)_{mn}] = \int_\tau \kappa(x; \mu) \nabla \psi_m^{\text{off}} \cdot \nabla \psi_n^{\text{off}} = R_{\text{off}}^T A(\mu) R_{\text{off}}
 \end{equation*}
\begin{equation*}
 \displaystyle S^{\text{on}}(\mu) = [s^{\text{on}}(\mu)_{mn}] = \int_\tau \widetilde{\kappa}(x; \mu) \psi_m^{\text{off}} \psi_n^{\text{off}} = R_{\text{off}}^T S(\mu) R_{\text{off}},
 \end{equation*}
and $\kappa(x; \mu)$ and $\widetilde{\kappa}(x; \mu)$ are now parameter dependent. To generate the online space we then choose the smallest $M_{\text{on}}$ eigenvalues from Eq.~\eqref{oneig} and form the corresponding eigenvectors in the offline space by setting
$\psi_k^{\text{on}} = \sum_j \Psi_{kj}^{\text{on}} \psi_j^{\text{off}}$ (for $k=1,\ldots, M_{\text{on}}$), where $\Psi_{kj}^{\text{on}}$ are the coordinates of the vector $\Psi_{k}^{\text{on}}$.

\subsection{Global coupling}
\label{globcoupling}

\subsubsection{Continuous Galerkin coupling}
\label{sec:CGcoupling}
In this subsection we aim to create an appropriate solution space and variational formulation that is suitable for a continuous Galerkin approximation of Eq.~\eqref{eq:original2}. We begin with an initial coarse space $V^{\text{init}}_0(\mu) = \text{span}\{ \chi_i \}_{i=1}^{N_v}$ (we use $N_v$ to denote the number of coarse vertices), where the $\chi_i$ are the standard multiscale partition of unity functions defined by
\begin{eqnarray} \label{pou}
-\text{div} \left( \kappa(x; \mu) \, \nabla \chi_i  \right) = 0 \quad K \in \omega_i \\
\chi_i = g_i \quad \text{on} \, \, \, \partial K, \nonumber
\end{eqnarray}
for all $K \in \omega_i$, where $g_i$ is assumed to be linear. Referring back to Eq.~\eqref{eqn:eigenmatrix} (for example), we note that the summed, pointwise energy $\widetilde{\kappa}$ required for the eigenvalue problems will be defined as
\begin{equation*}
\widetilde{\kappa} = \kappa \sum_{i=1}^{N_v} H^2 | \nabla \chi_i |^2.
\end{equation*}
We then multiply the partition of unity functions by the eigenfunctions in the online space $V_{\text{on}}^{\omega_i}$ to construct the resulting basis functions
\begin{equation} \label{cgbasis}
\psi_{i,k}^{\text{CG}} = \chi_i \psi_k^{\omega_i, \text{on}} \quad \text{for} \, \, \,
1 \leq i \leq N_v \, \, \,  \text{and} \, \, \, 1 \leq k \leq M_{\text{on}}^{\omega_i},
\end{equation}
where $M_{\text{on}}^{\omega_i}$ denotes the number of online eigenvectors that are chosen for each coarse node $i$. We note that the construction in Eq.~\eqref{cgbasis} yields inherently continuous basis functions due to the multiplication of online eigenvectors with the initial (continuous) partition of unity. This convention is not necessary for the discontinuous Galerkin global coupling, and is a focal point of contrast between the respective methods. However, with the continuous basis functions in place, we define the continuous Galerkin spectral multiscale space as
\begin{equation} \label{cgspace}
V_{\text{on}}^{\text{CG}}(\mu) = \text{span} \{ \psi_{i,k}^{\text{CG}}: \,  \, 1 \leq i \leq N_v \, \, \,  \text{and} \, \, \, 1 \leq k \leq M_{\text{on}}^{\omega_i}  \}.
\end{equation}
Using a single index notation, we may write $V_{\text{on}}^{\text{CG}}(\mu) = \text{span} \{ \psi_{i}^{\text{CG}} \}_{i=1}^{N_c}$, where $N_c$ denotes the total number of basis functions that are used in the coarse space construction. We also construct an operator matrix $R_0^T = \left[ \psi_1^{\text{CG}} , \ldots, \psi_{N_c}^{\text{CG}} \right]$ (where $\psi_i^{\text{CG}}$ are used to denote the nodal values of each basis function defined on the fine grid), for later use in this subsection.

Before introducing the continuous Galerkin formulation, we recall that the parameter $\mu$ is used to denote a solution that is computed at a previous iteration level (see Eq.~\eqref{eq:original2}). In turn, to update the solution at the current iteration level we seek $u^{\text{CG}}_{\text{ms}}(x; \mu) = \sum_i c_i \psi_i^{\text{CG}}(x; \mu) \in V_{\text{on}}^{\text{CG}}$ such that
\begin{equation} \label{cgvarform}
a^{\text{CG}}(u_{\text{ms}}^{\text{CG}}, v; \mu) = (f, v) \quad \text{for all} \,\,\, v \in V_{\text{on}}^{\text{CG}},
\end{equation}
where
$ \displaystyle a^{\text{CG}}(u, v; \mu) = \int_D \kappa(x;\mu) \nabla u \cdot \nabla v \, dx$, and $ \displaystyle (f,v) = \int_D f v \, dx$. We note that variational form in \eqref{cgvarform} yields the following linear algebraic system
\begin{equation}
A_0 U_0^{\text{CG}} = F_0,
\end{equation}
where $U^{\text{CG}}_0$ denotes the nodal values of the discrete CG solution, and
\begin{equation*}
A_0(\mu) = [a_{IJ}] = \int_D \kappa(x; \mu)  \, \nabla \psi_I^{\text{CG}} \cdot \nabla \psi_J^{\text{CG}} \, dx \quad \text{and} \quad F_0 = [f_I] = \int_D f \psi_I^{\text{CG}} \, dx.
\end{equation*}
Using the operator matrix $R_0^T$, we may write $A_0(\mu) = R_0 A(\mu) R_0^T$ and $F_0 = R_0 F$, where $A(\mu)$ and $F$ are the standard, fine scale stiffness matrix and forcing vector corresponding to the form in Eq.~\eqref{cgvarform}. We also note that the operator matrix may be analogously used in order to project coarse scale solutions onto the fine grid.

\subsubsection{Discontinuous Galerkin coupling}
\label{dgcoupling}
One can also use the discontinuous Galerkin (DG) approach
(see also \cite{MR1885715,MR2002258,MR2431403})
to couple multiscale basis functions. This may avoid the use of
the partition of unity functions;
however, a global formulation needs to be chosen carefully.
We have been investigating the use of DG coupling and the detailed
results will be presented elsewhere, see \cite{DGcoupling}. 
Here, we would like to
briefly mention a general global coupling that can be used.
The global formulation is given by
\begin{equation}\label{eq:disc}
{a}^{\text{DG}}(u, v;\mu) = f(v) \quad \mbox{ for all }
\quad v=
\{v_{\text{K}} \in \text{V}^\text{K}\} ,
\end{equation}
where
\begin{equation}\label{eq:def:a-h}
{a}^{\text{DG}}(u, v;\mu) =\sum_{\text{K}} {a}^{\text{DG}}_\text{K}(u,v;\mu)~~~
\mbox{and}~~~~f(v) = \sum_{\text{K}} \int_{\text{K}} f v_\text{K} dx,
\end{equation}
for all
$u=\{u_\text{K}\}, v=\{v_\text{K}\}$ with \text{K} being the coarse element depicted in Figure \ref{schematic}.
Each local bilinear form ${a}^{\text{DG}}_\text{K}$ is given as a sum of three
bilinear forms:
\begin{equation}\label{eq:def:a^hat-i}
 {a}^{\text{DG}}_K(u,v;\mu) :=  a_K(u, v;\mu) + r_K(u, v;\mu) + p_K(u, v;\mu),
\end{equation}
where $a_\text{K}$ is the bilinear form,
\begin{equation}\label{eq:def:a-i}
a_K(u, v;\mu) := \int_{K}\!\!\kappa_\text{K}(x;\mu) \nabla u_\text{K} \cdot \nabla v_\text{K} dx,
\end{equation}
where $\kappa_K(x;\mu)$ is the restriction of $\kappa(x;\mu)$ in $K$;
the $r_\text{K}$ is the symmetric bilinear form,
\begin{equation}\label{eq:def:s-i}
r_\text{K}(u, v;\mu) := \sum_ {E\subset \partial K}
 \frac{1}{l_{E}}\int_{E}  \widetilde{\kappa}_{E}(x;\mu) \left( \frac{\partial
  u_\text{K}}{\partial n_\text{K}} (v_\text{K} - v_{\text{K}'}) + \frac{\partial v_\text{K}}{\partial n_\text{K}} (u_{\text{K}'} - u_{\text{K}})
  \right) ds,
\nonumber
\end{equation}
where $ \widetilde{\kappa}_{E}(x;\mu)$ is the harmonic average of $\kappa(x;\mu)$
along the edge $E$, $l_E=1$ if E is on the boundary of the macrodomain, and $l_E=2$ if
$E$ is
an inner edge of the macrodomain.
Here, $\text{K}'$ and $\text{K}$ are two coarse-grid elements
sharing the common edge $E$;
and $p_K$ is the penalty bilinear form,
\begin{equation} \label{eq:def:p-i}
p_\text{K}(u, v;\mu) := \sum_{E \subset \partial K}
\frac{1}{l_E}\frac{1}{h_E} {\delta_E}
\int_{E}
\widetilde{\kappa}_{E}(x;\mu) (u_{K'} - u_K)(v_{K'} - v_K)ds.
\end{equation}
Here  $h_E$ is harmonic average of the length of the edge $E$ and $E'$,
$\delta_E$ is a positive penalty parameter that needs to be selected
and its choice affects the performance of GMsFEM.
One can choose other eigenvalue problems within the DG framework. 
See \cite{DGcoupling}.

As mentioned before that for Discontinuous Galerkin formulation, the
support of basis functions are coarse element K as depicted
in Figure \ref{schematic}. Besides, the inherent unconformal property of DG formulation
 determines the removal of the partition of unity functions while constructing basis functions in Equation \eqref{cgbasis}. Similarly, we can obtain the discontinuous Galerkin spectral multiscale space as
\begin{equation} \label{dgspace}
V_{\text{on}}^{\text{DG}}(\mu) = \text{span} \{ \psi_{k}^{\text{DG}}: \,  \, \, 1 \leq k \leq M_{\text{on}}^{K}  \},
\end{equation}
For every coarse element K.

Using the same process in the continuous Galerkin formulation, we can obtain an operator matrix constructed by the basis functions of $V_{\text{on}}^{\text{DG}}(\mu)$. For the consistency of the notation, we denote the matrix as $R_{0}$, and $R_0^T = \left[ \psi_1^{\text{DG}} , \ldots, \psi_{N_c}^{\text{DG}} \right]$. Recall that $N_c$ denote the total number of coarse basis functions.

Solving the problem \eqref{eq:original} in the coarse space $V_{\text{on}}^{\text{DG}}(\mu)$ using the DG formulation described in Equation \eqref{eq:disc} is equivalent to seeking
$u^{\text{DG}}_{\text{ms}}(x; \mu) = \sum_i c_i \psi_i^{\text{DG}}(x; \mu) \in V_{\text{on}}^{\text{DG}}$ such that
\begin{equation} \label{dgvarform}
a^{\text{DG}}(u_{\text{ms}}^{\text{DG}}, v; \mu) = f( v) \quad \text{for all} \,\,\, v \in V_{\text{on}}^{\text{DG}},
\end{equation}
where
$ \displaystyle a^{\text{DG}}(u, v; \mu) $ and $f( v)$ are defined in Equation \eqref{eq:def:a-h}.
Similar as the CG case, we can obtain a coarse linear algebra system
\begin{equation}
A_0 U_0^{\text{DG}} = F_0,
\end{equation}
where $U^{\text{DG}}_0$ denotes the discrete coarse DG solution, and
\begin{equation*}
A_0(\mu) = R_0 A(\mu) R_0^T \quad \text{and} \quad F_0 = R_0 F,
\end{equation*}
where $A(\mu)$ and $F$ are the standard, fine scale stiffness matrix and forcing vector corresponding to the form in Eq.~\eqref{eq:def:a-h}. After solving the coarse matrix, we can use the operator matrix $R_0$ to obtain the fine-scale solution in the form of $R_0^{T}U_0^{\text{DG}}$.

\section{Numerical Results}
\label{numerical}

\begin{figure}\centering
 \subfigure[Random field values ]{\label{fig:perm1}
    \includegraphics[width = 0.45\textwidth, keepaspectratio = true]{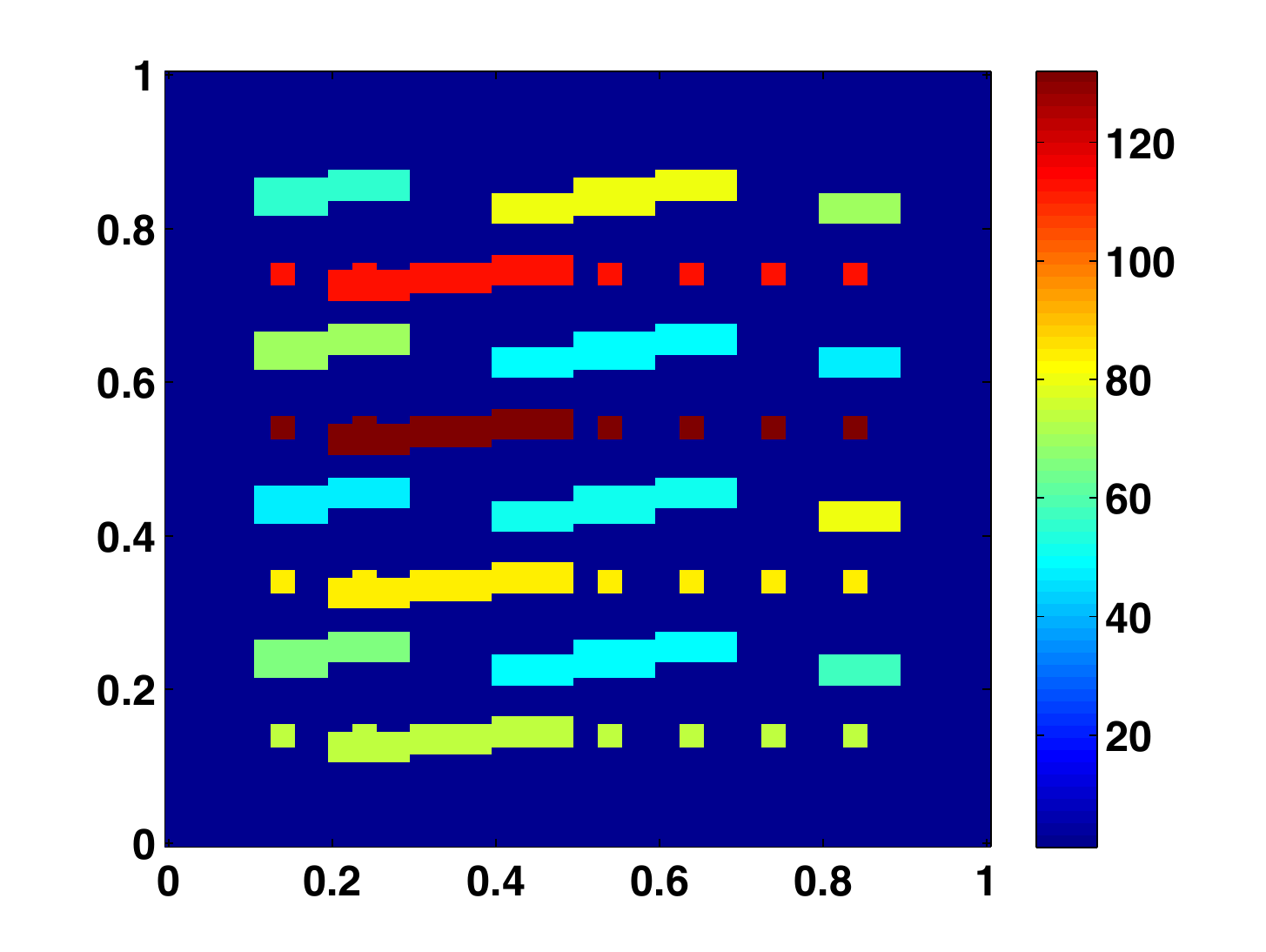}
   }
  \subfigure[Fixed field values]{\label{fig:perm2}
     \includegraphics[width = 0.45\textwidth, keepaspectratio = true]{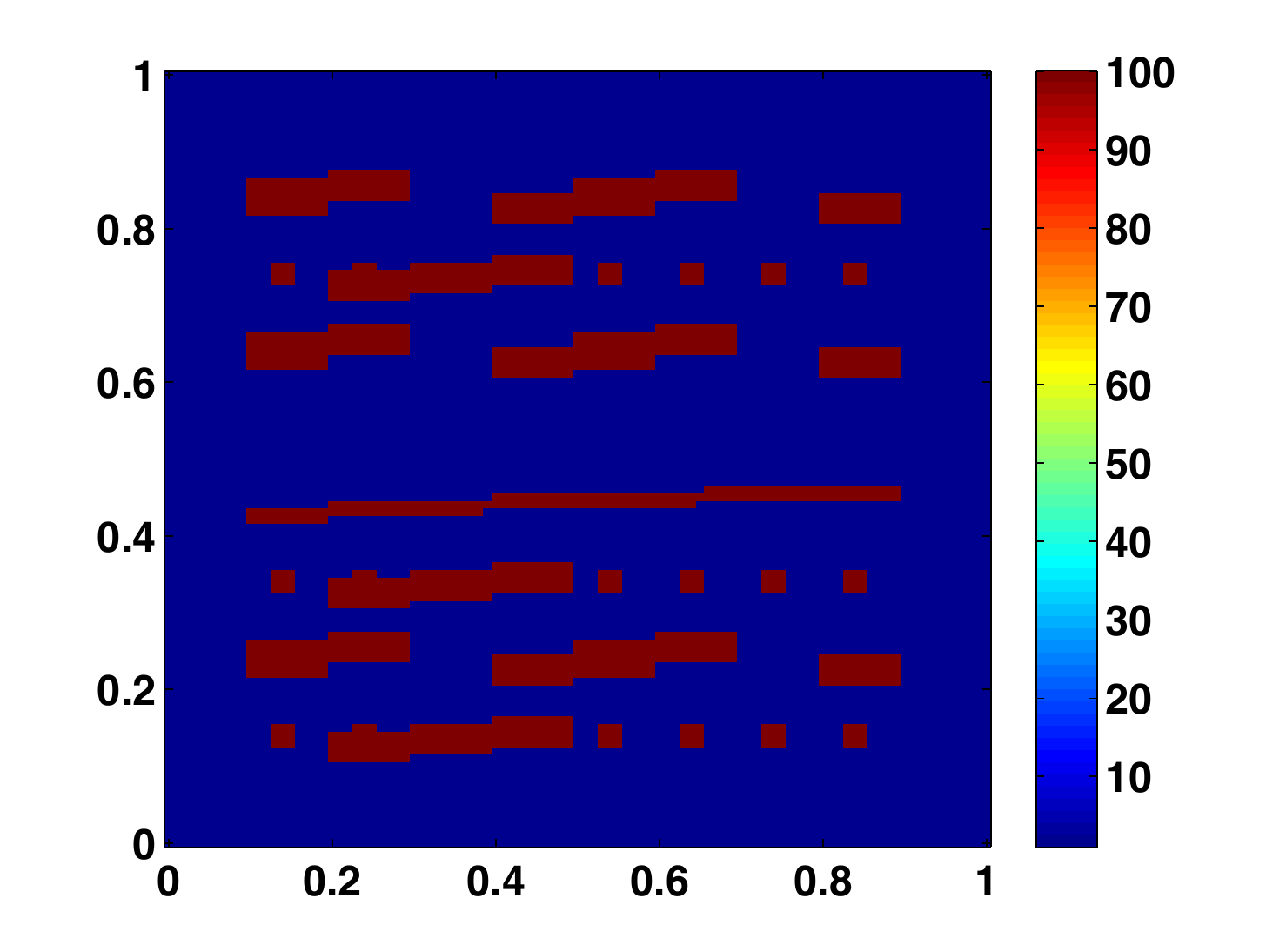}
  }
 \caption{High-contrast permeability fields}\label{fig:perm}
\end{figure}

In this section we solve the nonlinear, elliptic model equation given in Eq.~\eqref{eq:original1} using both the continuous (CG) and discontinuous Galerkin (DG) GMsFEM formulations described in Sect.~\ref{cgdgmsfem}. More specifically, we consider the equation
\begin{subequations}\label{eqn:problem}
\begin{align}
 -\DIV\big(\mathrm{e}^{\kappa(x)u(x)}\nabla u(x)\big)&=f \, \, \, \text{in}  \, \,   D\\
u&=0 \, \, \, \text{on} \;\partial D,
\end{align}
\end{subequations}
where the general coefficient from \eqref{eq:original1} is taken to be  $\kappa(x; u) = \mathrm{e}^{\kappa(x)u(x)}$.
For the coefficient $\kappa(x)$, we consider the high-contrast permeability fields as illustrated Fig.~\ref{fig:perm}.  Fig.~\ref{fig:perm1} represents a field whose high-permeability values are randomly assigned, while the field in Fig.~\ref{fig:perm2} has a different channelized structure with fixed maximum values. We use a source term $f=0.1$, and solve the problem on the unit two-dimensional domain $D = [0,1]{\times}[0,1]$.

To solve Eq.~\eqref{eqn:problem} we first linearize it by using a Picard iteration.
In particular, for a given initial guess
$u^{0}$ we solve
\begin{subequations}
\begin{align}
-\DIV\big(\mathrm{e}^{\kappa(x)u^{n}(x)}\nabla u^{n+1}(x) \big)&=f \, \, \, \text{in} \; D\\
u^{n+1}&=0 \, \, \,  \text{on} \;\partial D,
\end{align}
\end{subequations}
for $n\geq 0$.

In our simulations, we take the initial guess $u^{0}=0$, and terminate the iterative loop
when
$\norm{A(u^{n+1})u^{n+1}-b}\leq \delta \norm{b}$, where $\delta$ is the tolerance
for the iteration and we select $\delta=10^{-3}$. We note that $A$ and $b$ correspond to
the linear system resulting from either the CG or DG global formulations. In particular,  we solve the problem as follows:
\begin{align}\label{eqn:linearsystem}
A(u^{n})u^{n+1}=b \, \, \, \text{for} \, \, n=0,\,1,\ldots.
\end{align}
We note that since $u^{n}$ and $u^{n+1}$ will not necessarily be computed in coarse spaces of the same dimension, we cannot directly use the residual criterion
listed above. Actually, we use the Galerkin projection of the fine
solution to the corresponding coarse space to calculate the residual
error from above.

\begin{figure}\centering
 \subfigure[$\kappa_1(x)$]{\label{fig:permi}
    \includegraphics[width = 0.45\textwidth, keepaspectratio = true]{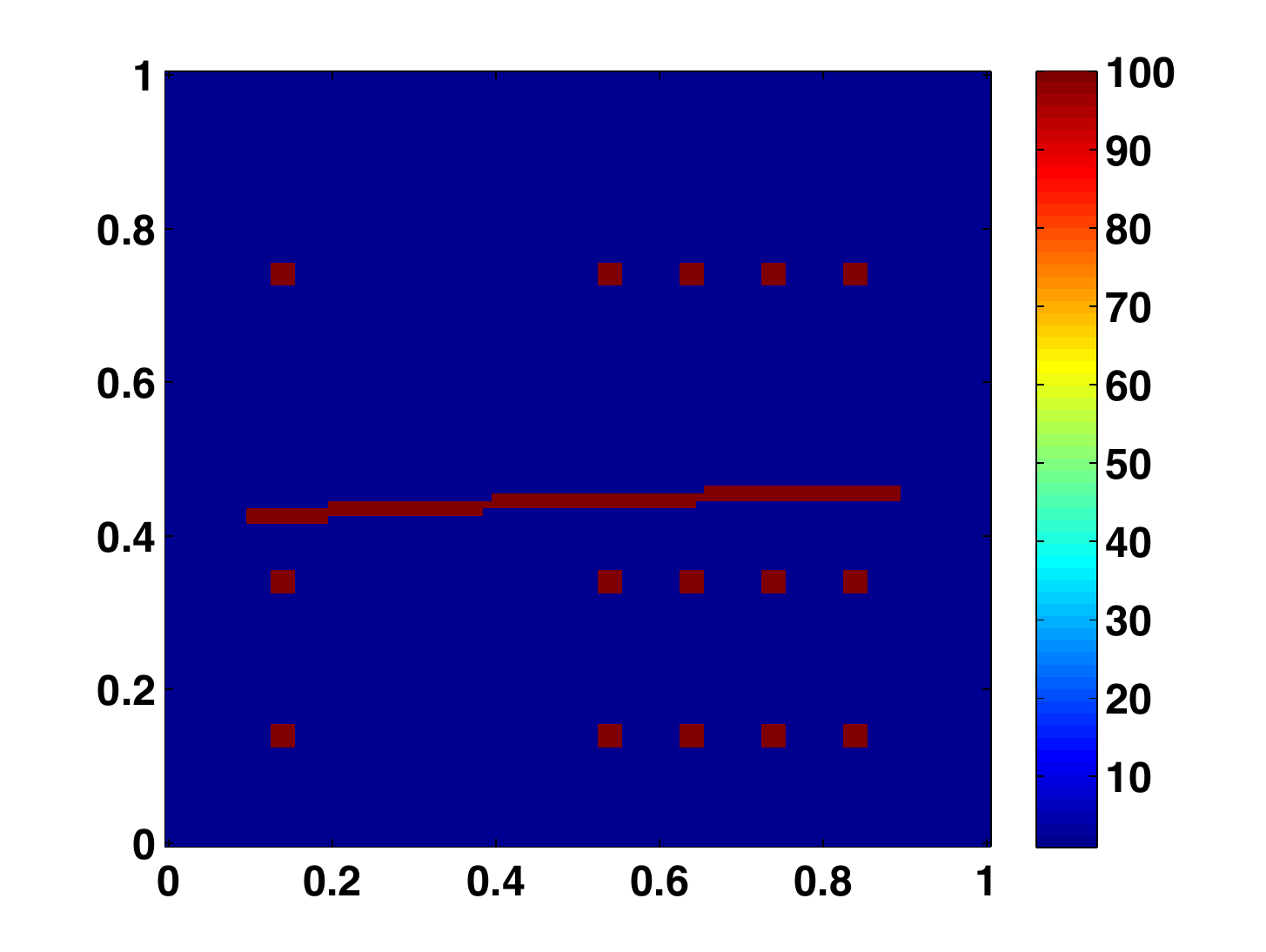}
   }
  \subfigure[$\kappa_2(x)$]{\label{fig:permii}
     \includegraphics[width = 0.45\textwidth, keepaspectratio = true]{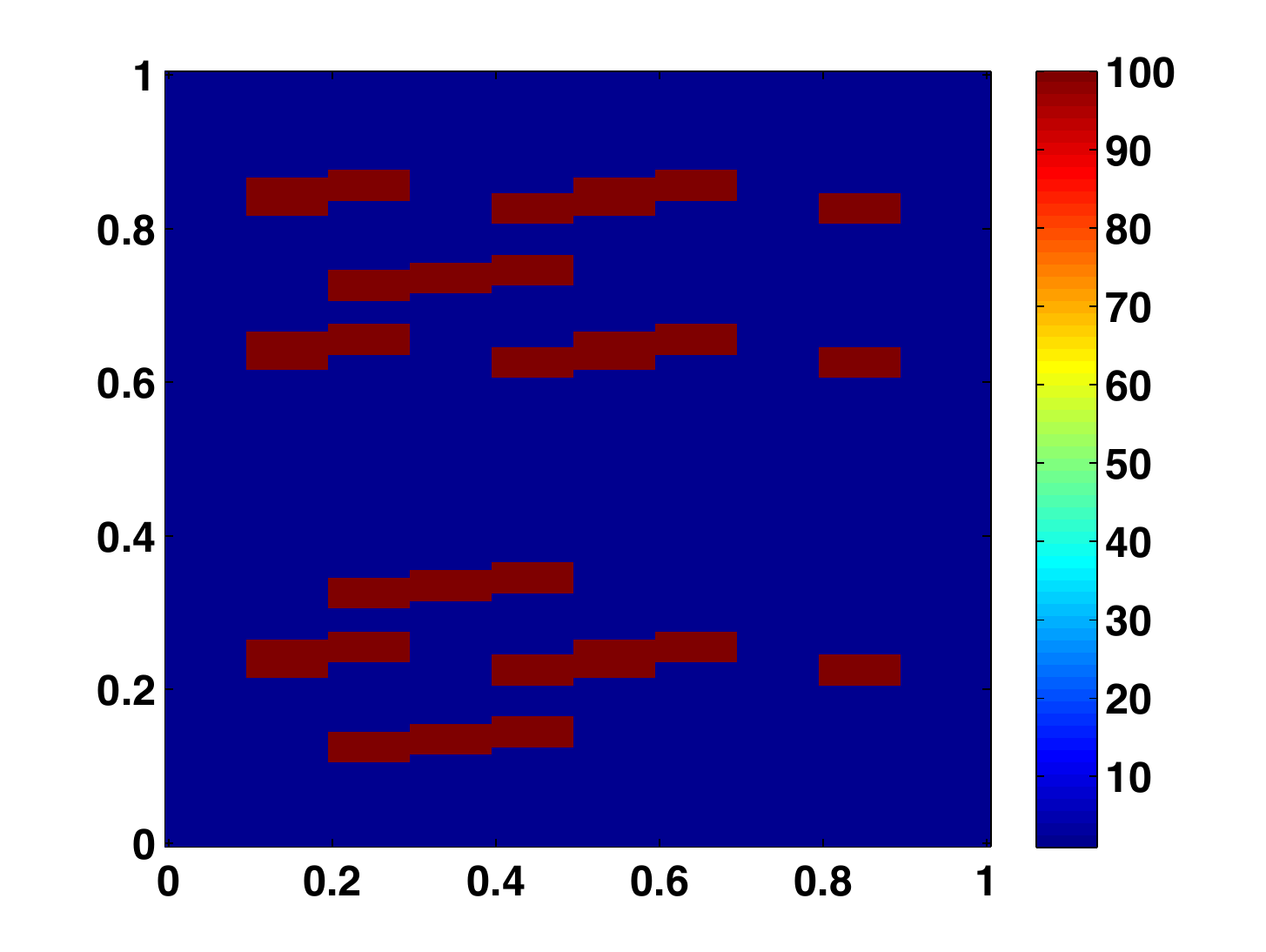}
  }
 \caption{Decomposition of permeability field \ref{fig:perm2}}\label{fig:decofperm}
\end{figure}

\begin{remark}
In this section we will consider two types of coefficients $\kappa(x)$ to be used in  Eq.~\eqref{eqn:problem}. We recall that throughout the paper we have used an auxiliary variable $\mu = \overline{u}^n$ to denote the solution dependence of the nonlinear problem. As such, we have referred to the model equation as parameter-dependent while describing the iterative solution procedure. Consequently, we are careful to introduce (and distinguish) a related case where we use a ``physical" parameter $\mu^{\text{p}}$ for the purpose of constructing a field of the form $\kappa(x) = \mu^{\text{p}} \kappa_1(x) + (1-\mu^{\text{p}}) \kappa_2(x)$. See Fig.~\ref{fig:decofperm} for an illustration of $\kappa_1(x)$ and $\kappa_2(x)$. We note that the coefficient will be constructed by summing contributions that depend on the physical parameter $\mu^{\text{p}}$, in addition to the auxiliary parameter dependence from the iterative form. In Subsect. \ref{noparameter} we use a field that does not depend on $\mu^{\text{p}}$, and in Subsect. \ref{withparameter} we use a field that does depend on $\mu^{\text{p}}$.
\end{remark}

\subsection{Parameter-independent permeability field}
\label{noparameter}
In the following simulations we
first generate a snapshot space, use a spectral decomposition to
obtain the offline space, and then for an initial guess apply a similar spectral decomposition to obtain the online space. We recall that in order to construct the snapshot space we choose a specified number of eigenfunctions (denoted by $M_{\text{snap}}$) on either a coarse neighborhood or coarse element depending on whether we use continuous (CG) or discontinuous Galerkin (DG) global coupling, respectively. In our simulations, we select the range of solutions $[u_\text{min}, u_\text{max}]$ that correspond to solving the fine scale equation using a source term that ranges from $f \in [0.1, 1]$. For the first set of simulations we divide the domain $[u_\text{min}, u_\text{max}]$ into $N_s-1$ equally spaced subdomains to obtain $N_s$ discrete points
$u_1, \ldots, u_{N_s}$. For these simulations we fix a value of $N_s=9$.

For either formulation, we solve a localized eigenvalue problems as defined in Subsect. \ref{locbasis} for each point $u_j$ on a coarse neighborhood and keep a specified number of eigenfunctions. For example, in the CG case we keep $l_{\text{max}} = 3$ snapshot eigenfunctions, and this construction leads to a local space of dimension $M_{\text{snap}} = l_{\text{max}} {\times} N_s = 3 {\times} 9 = 27$. In the DG case, we adaptively choose the number of eigenfunctions based on a consideration of the eigenvalue differences. In the offline space construction we fix $\overline{u}$ as the average of the previously defined fixed snapshot values. We then solve the offline eigenvalue problem and construct the offline space by keeping the eigenvectors corresponding to a specified number of dominant eigenvalues. At the online stage we use the initial guess $u^0=0$ in order to solve the respective eigenvalue problem required for the space construction. We note that the size of our online space and the associated solution accuracy will depend on the number of eigenvectors that we keep in the online space construction.

In the CG formulation, we recall that the online eigenfunctions are multiplied by the corresponding partition of unity functions with support in the same neighborhood of the respective coarse node. We then solve Eq.~\eqref{eqn:problem} iteratively within the online space. In particular, for each iteration we update the online space and
solve the equation Eq.~\eqref{eqn:problem} using the previously computed solution.

\begin{figure}[htb]\centering

    \includegraphics[width = 0.33\textwidth, keepaspectratio = true]{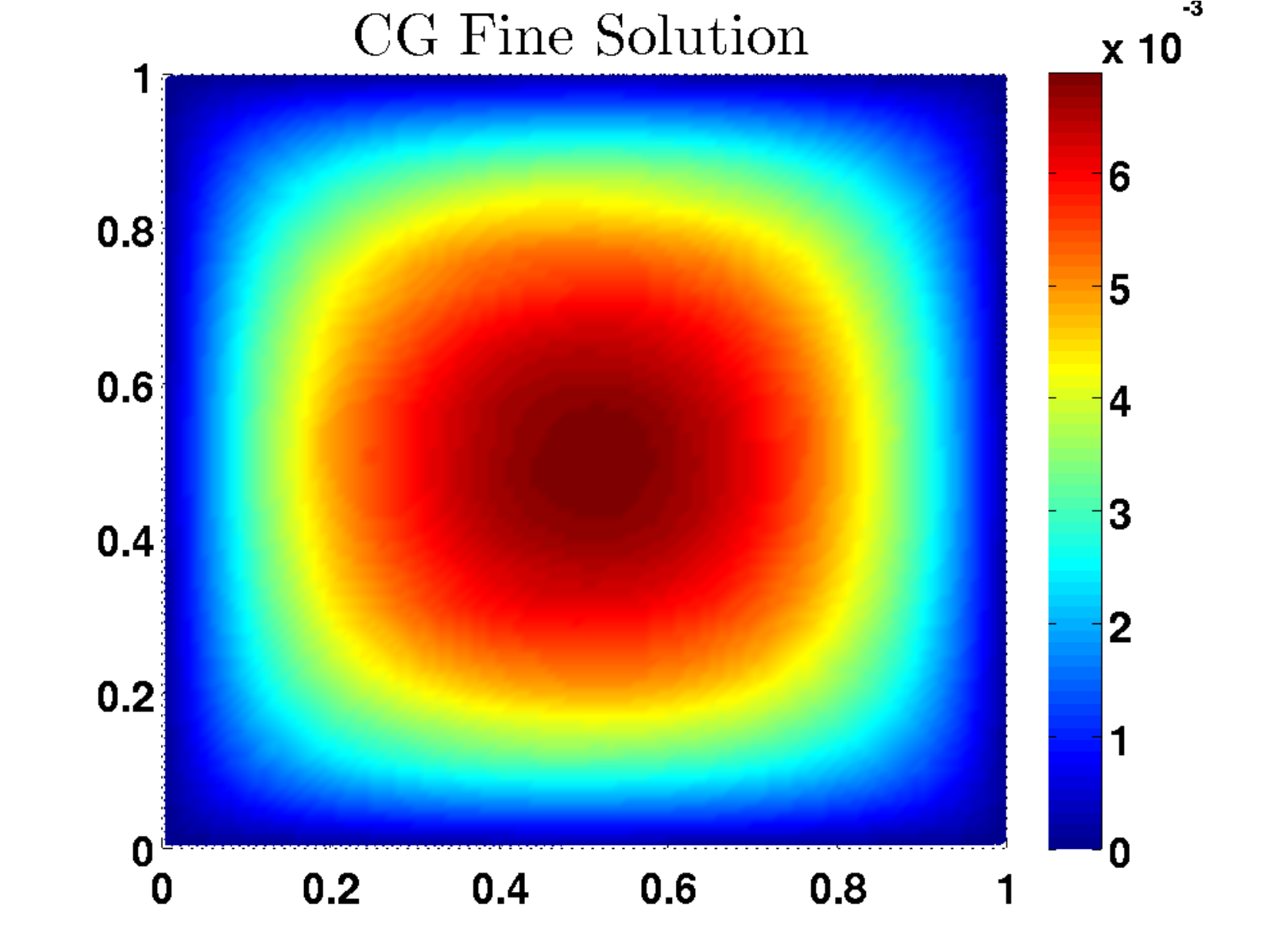}
    \hspace*{-0.5cm} 
    \includegraphics[width = 0.33\textwidth, keepaspectratio = true]{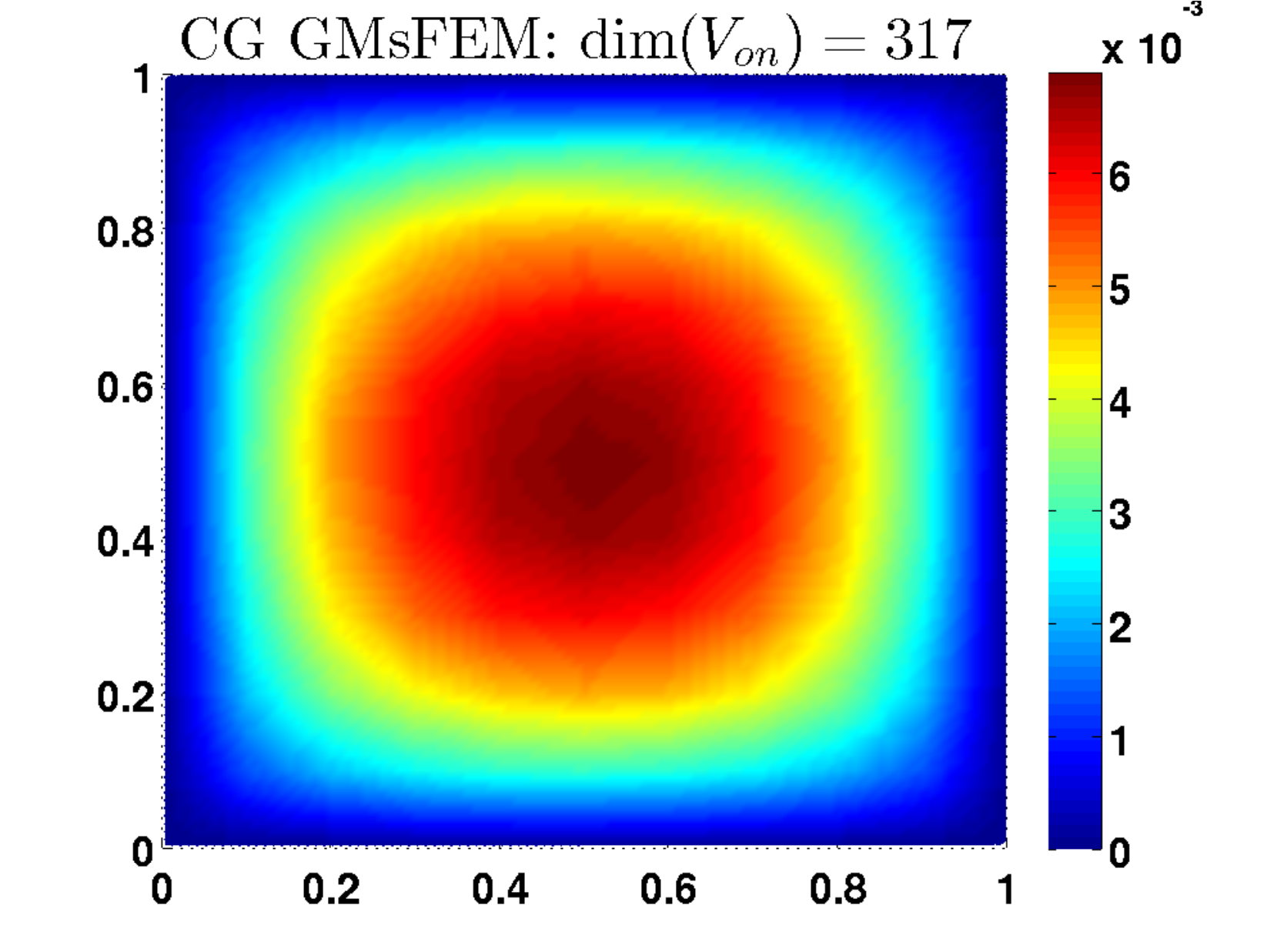}
    \hspace*{-0.5cm} 
    \includegraphics[width = 0.33\textwidth, keepaspectratio = true]{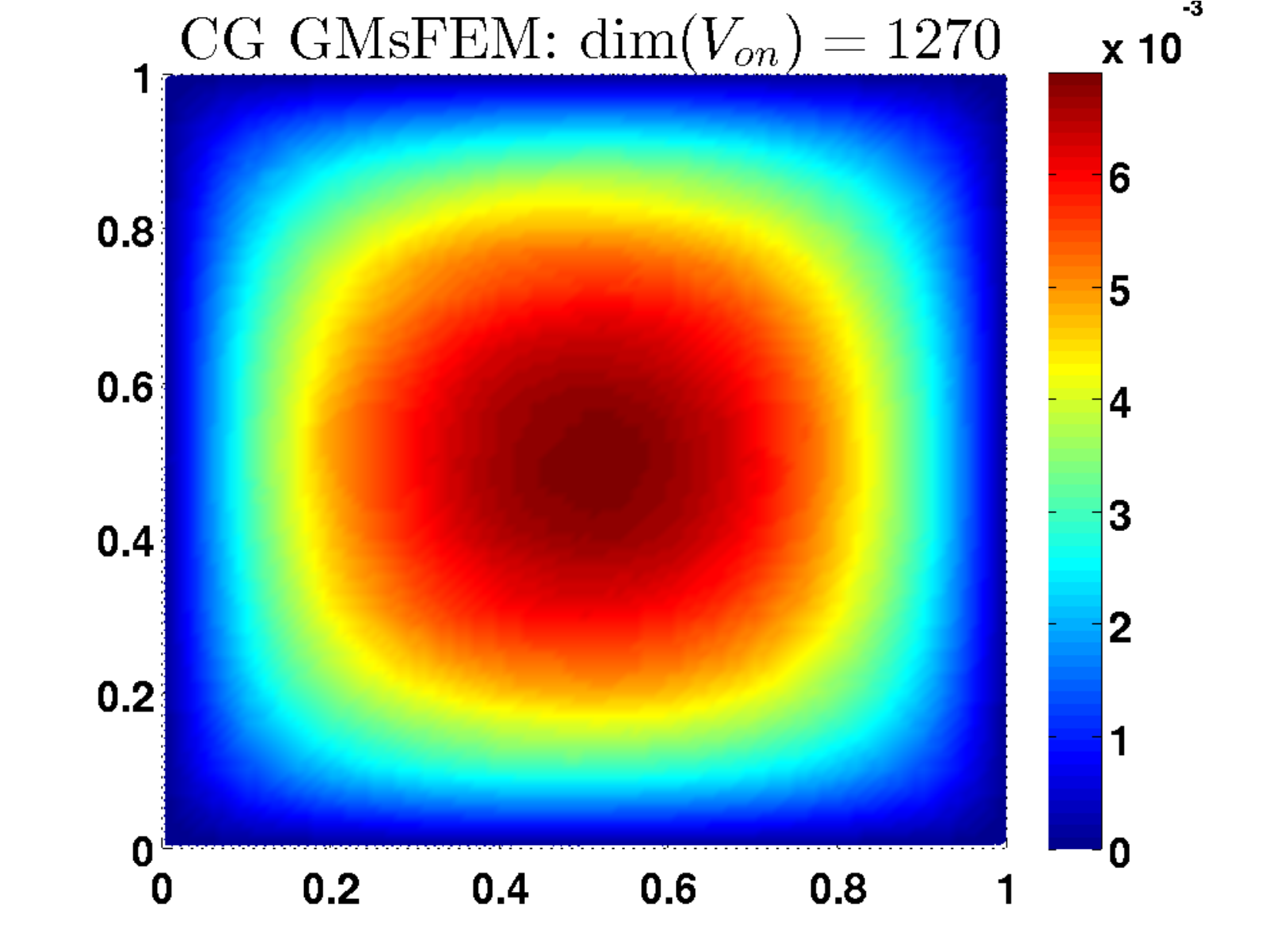}

 \caption{Comparison of fine and coarse CG solutions correpsonding to
               Fig.~\ref{fig:perm2}}
 \label{fig:cgsolnperm2}
\end{figure}

In the simulations using the CG formulation we discretize our domain into coarse elements of size $H=1/10$, and fine elements of size $h=1/100$. The results corresponding to the  permeability fields from Figs.~\ref{fig:perm1} and \ref{fig:perm2} are shown in Tables
\ref{table:cgperm1} and \ref{table:cgperm2}, respectively. The first column shows the dimension of the online solution space, and the second
column shows the eigenvalue $\lambda^*$ which corresponds to the first eigenfunction that is discarded from space enrichment. We note that this eigenvalue is an important consideration in error estimates of enriched multiscale spaces (\cite{egw10}). As a formal consideration, we mention that the error analysis typically yields estimates of the form
$\| u - u_{\text{ms}} \| \sim \mathcal{O}(H^{\gamma} \lambda^*)$ when the dominant eigenvalues are taken to be small. The next two
columns correspond to the $L^{2}$-weighted relative error
$\norm{u-u_{\text{ms}}}_{L^2_\kappa(D)} / \norm{u}_{L^2_\kappa(D)} {\times} 100\%$ and energy relative error $\norm{u-u_{\text{ms}}}_{H^1_\kappa(D)} / \norm{u}_{H^1_\kappa(D)} {\times} 100\%$ between the GMsFEM solution $u_{\text{ms}}$ and the fine-scale solution $u$. We note that as the dimension of the online space increases (i.e., we keep more eigenfunctions in the space construction), the relative errors decrease accordingly. As an example, for the field in Fig.~\ref{fig:perm1}, we encounter $L^2$ relative errors that decrease from $1.43 - 0.24\%$, and energy relative errors that decrease from $16.12 - 6.85 \%$ as the online space is systematically enriched. In the tables, analogous errors between the online GMsFEM solution and the offline solution are computed. The dimension of the offline space is taken to be the maximum dimension of the online space. We note that in this case the Picard iteration converges in $4$ steps for all simulations. In Fig.~\ref{fig:cgsolnperm2} we also plot the fine and coarse-scale CG solutions that correspond to the field in Fig.~\ref{fig:perm2}. We note that the fine solution, and the coarse solutions corresponding to the largest and smallest online spaces are nearly indistinguishable.

\begin{figure}\centering
 \subfigure[Corresponds to Fig.~\ref{fig:perm1}]{\label{fig:CG21}
    \includegraphics[width = 0.45\textwidth, keepaspectratio = true]{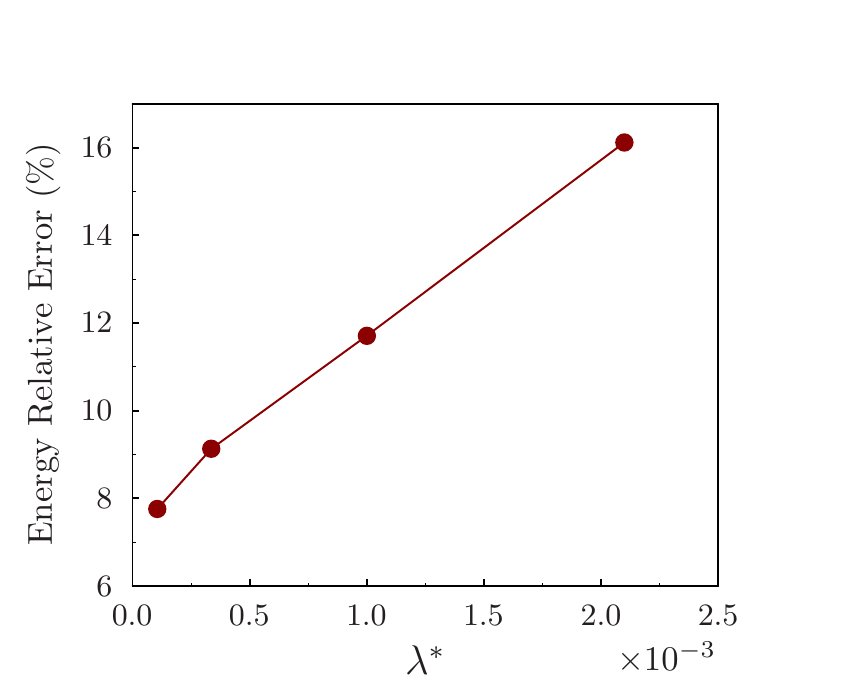}
   }
  \subfigure[Corresponds to Fig.~\ref{fig:perm2}]{\label{fig:CG22}
     \includegraphics[width = 0.45\textwidth, keepaspectratio = true]{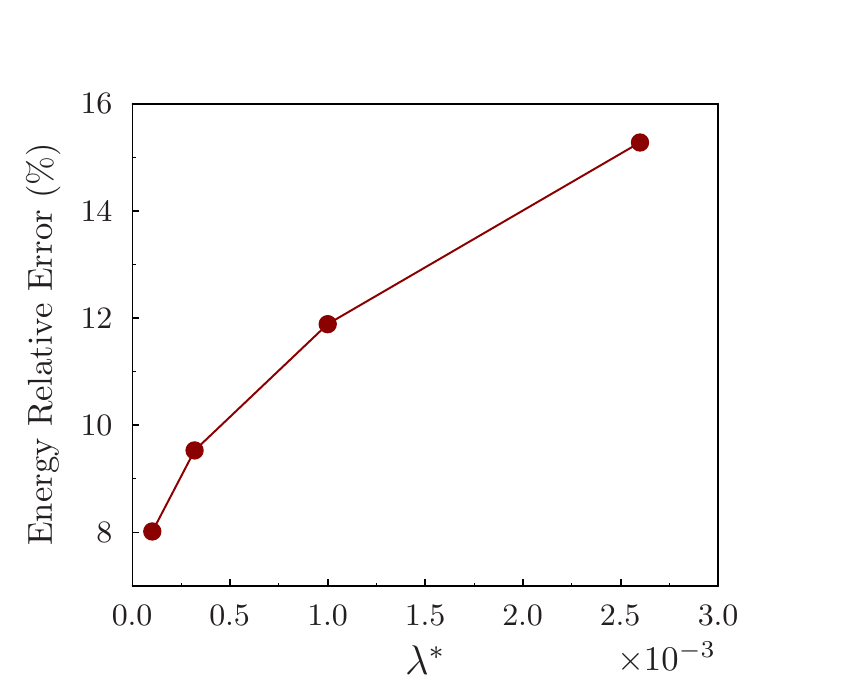}
  }

 \caption{Relation between the first discarded eigenvalue and the CG relative energy error;
               permeability from Fig.~\ref{fig:perm1} (left), permeability from Fig.~\ref{fig:perm2} (right) }
  \label{fig:CGorder2}
\end{figure}

We also illustrate the relation between the energy relative errors and
$\lambda^*$ in Fig.~\ref{fig:CGorder2} for the same permeability
fields considered above. From the plots in Fig.~\ref{fig:CGorder2}, we see
that the energy relative error predictably decreases as $\lambda^*$ decreases, thus following the appropriate error behavior.

\begin{table}
\tabcolsep 0pt %
\vspace*{1.5pt}
\begin{center}
\small \addtolength{\tabcolsep}{8.9pt}
\def\temptablewidth{1.0\textwidth}
\begin{tabular*}{\temptablewidth}{|c|c|c|c|c|c|}
\hline \multirow{2}{*}{$\dim(V_{\text{on}}^{\text{CG}})$}
&\multirow{2}{*}{$\lambda^{*}$}&
\multicolumn{2}{c|}{\;\;\;GMsFEM Relative Error (\%)\;\;\;}&\multicolumn{2}{c|}{\hspace{0cm}Online-Offline Relative Error (\%)\hspace{-.5cm}}\\
\cline{3-6} {}&{}&\hspace{.6cm}$L^{2}_{\kappa}(D)\hspace{.6cm}$ &
$H^{1}_{\kappa}(D)$&$\hspace{.6cm}L^{2}_{\kappa}(D)\hspace{.6cm}$&
$H^{1}_{\kappa}(D)$
\\
\hline\hline
       $319$&$0.0021$ &$1.43$&$16.12$ &$1.25$&$16.33$ \\

\hline
      $497$&$0.0010$ & $0.69$&$11.71$& $0.48$&$10.66$ \\
\hline
       $770$&$3.36\times 10^{-4}$& $0.40$&$9.13$ & $0.20$&$7.30$\\
\hline
      $1043$&$1.06\times 10^{-4}$&$0.31$&$7.76$ &$0.09$&$4.43$ \\
\hline
       $1270$& ---&$0.24$&$6.85$&0.00&0.00\\
\hline\hline
\end{tabular*}
      \end{center}
\caption{CG relative errors corresponding to the permeability field in
              Fig.~\ref{fig:perm1} }
\label{table:cgperm1}
       \end{table}
%
\begin{table}
\tabcolsep 0pt %
\vspace*{1.5pt}
\begin{center}
\small \addtolength{\tabcolsep}{8.9pt}
\def\temptablewidth{1\textwidth}
\begin{tabular*}{\temptablewidth}{|c|c|c|c|c|c|}
\hline \multirow{2}{*}{$\dim(V_{\text{on}}^{\text{CG}})$}
&\multirow{2}{*}{$\lambda^{*}$}&
\multicolumn{2}{c|}{\;\;\;GMsFEM Relative Error (\%)\;\;\;}&\multicolumn{2}{c|}{\hspace{0cm}Online-Offline Relative Error (\%)\hspace{-.5cm}}\\
\cline{3-6} {}&{}&\hspace{.6cm}$L^{2}_{\kappa}(D)\hspace{.6cm}$ &
$H^{1}_{\kappa}(D)$&$\hspace{.6cm}L^{2}_{\kappa}(D)\hspace{.6cm}$&
$H^{1}_{\kappa}(D)$
\\
\hline\hline
       $316$&$0.0026$ &$1.36$&$15.28$ &$1.18$&$15.74$ \\

\hline
      $482$&$0.0010$ & $0.71$&$11.89$& $0.51$&$11.17$ \\
\hline
       $722$&$3.18\times 10^{-4}$ & $0.43$&$9.53$ & $0.22$&$7.77$\\
\hline
      $996$&$1.02\times 10^{-4}$&$0.33$&$8.02$ &$0.11$&$4.72$ \\
\hline
       $1236$&---&$0.26$&$7.05$&0.00&0.00\\
\hline\hline
\end{tabular*}
       \end{center}
\caption{CG relative errors corresponding to the permeability field in
              Fig.~\ref{fig:perm2} }
\label{table:cgperm2}
\end{table}

In order to solve the model problem using the DG formulation, we note that the space of snapshots is constructed in a slightly different fashion. In this case, the selection of eigenvectors hinges on a comparison between the difference of consecutive eigenvalues resulting from the localized computations. In contrast to the CG case, the initial number of eigenfunctions (call this number $l_{\text{init}}^K$) used in the snapshot space construction are adaptively chosen based on the relative size of consecutive eigenvalues. For the results corresponding to the DG formulation, we note that two configurations for the snapshot space construction are used. In particular, we consider a case when the original number of eigenfuctions $l_{\text{init}}^K$ are used in the construction, and a case when $l_{\text{max}}^K = l_{\text{init}}^K + 3$ are used in the construction.

\begin{figure}[htb]\centering

  \includegraphics[width = 0.33\textwidth, keepaspectratio = true]{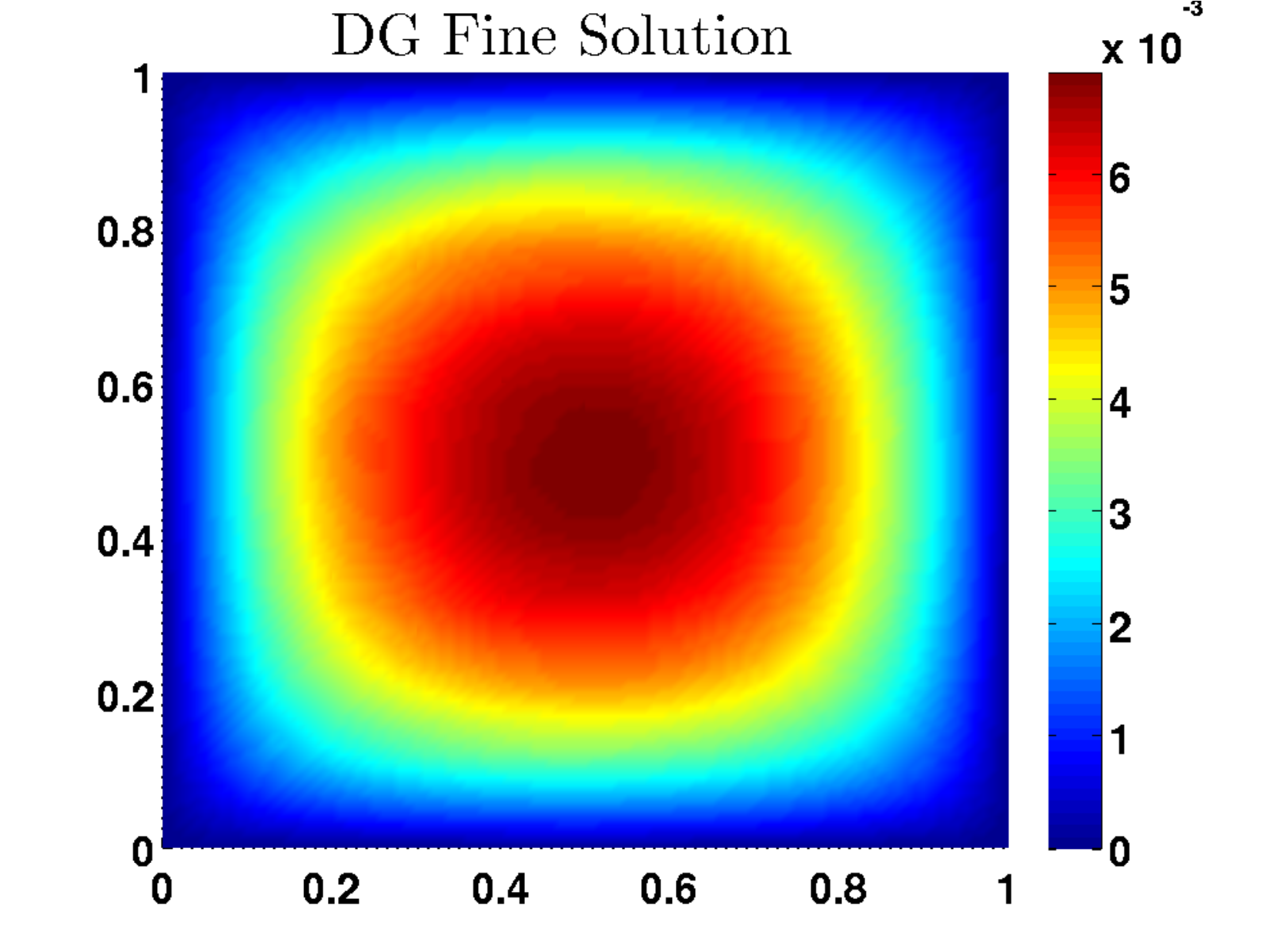}
    \hspace*{-0.5cm} \includegraphics[width = 0.33\textwidth, keepaspectratio = true]{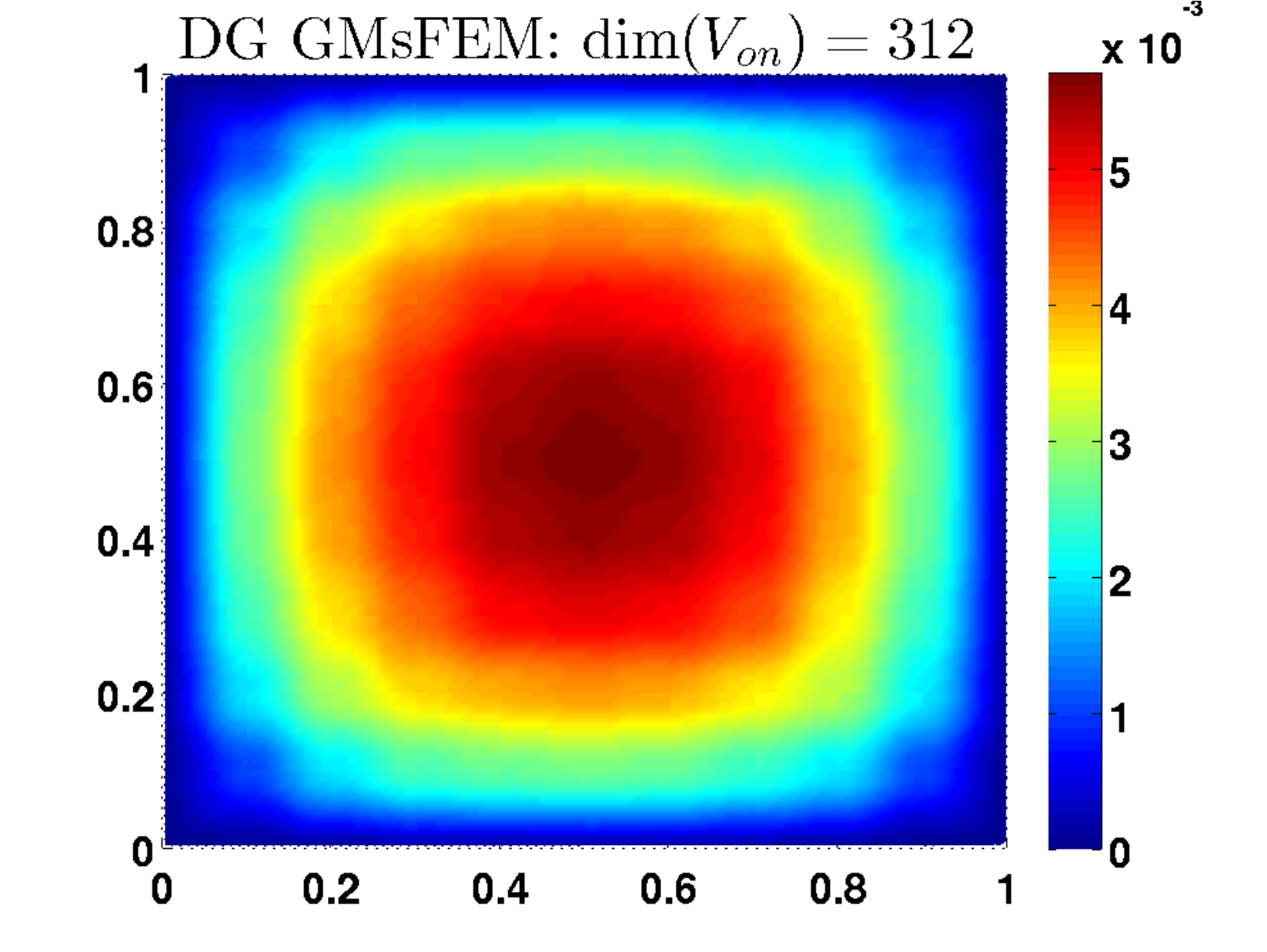}
    \hspace*{-0.5cm} \includegraphics[width = 0.33\textwidth, keepaspectratio = true]{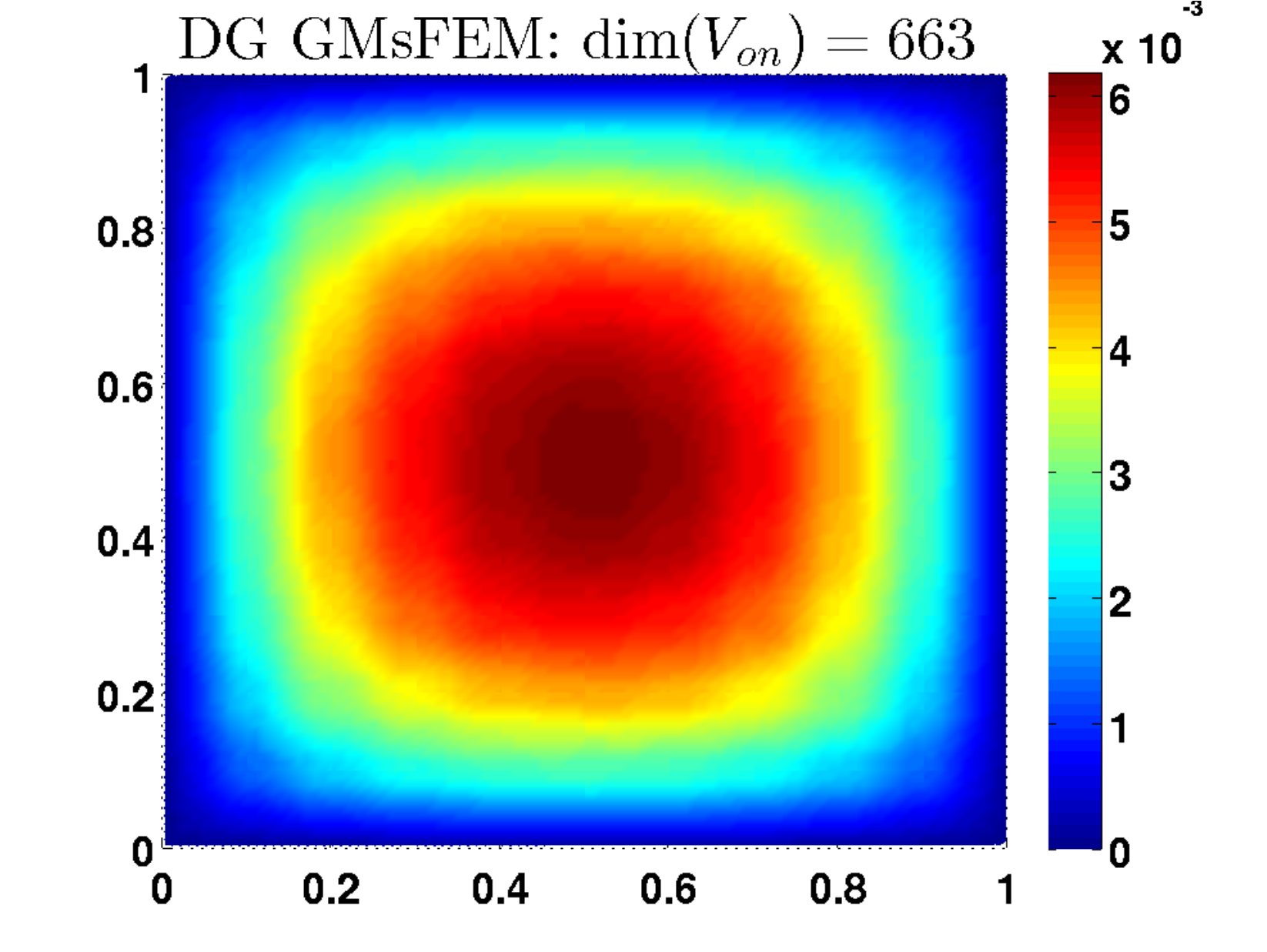}

 \caption{Comparison of fine and coarse DG solutions correpsonding to
               Fig.~\ref{fig:perm2}}

 \label{fig:dgsolnperm2}
\end{figure}

\begin{figure}\centering
 \subfigure[Corresponds to Fig.~\ref{fig:perm1}]{\label{fig:DG2I1}
    \includegraphics[width = 0.45\textwidth, keepaspectratio = true]{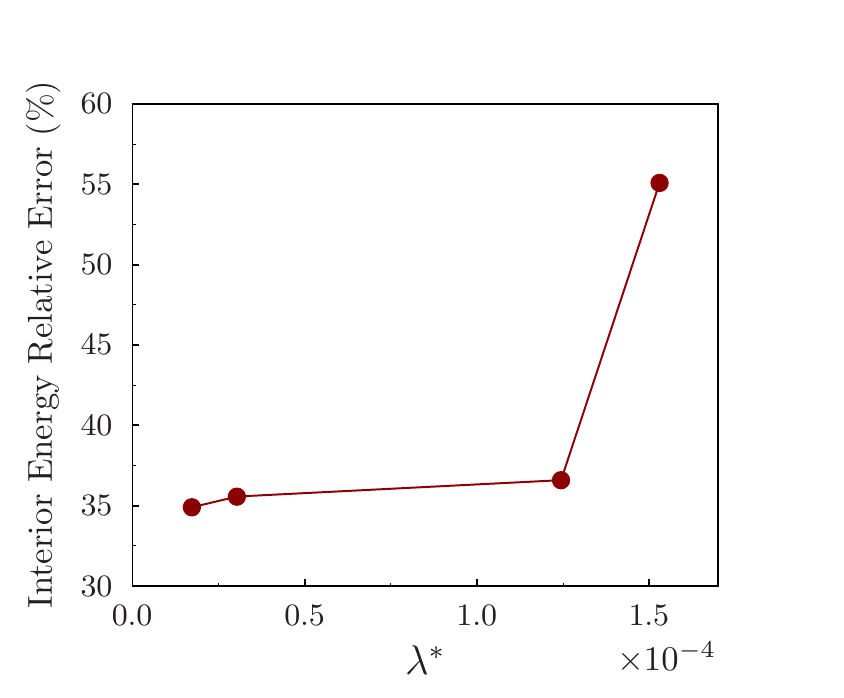}
   }
\subfigure[Corresponds to Fig.~\ref{fig:perm2}]{\label{fig:DG2I2}
    \includegraphics[width = 0.45\textwidth, keepaspectratio = true]{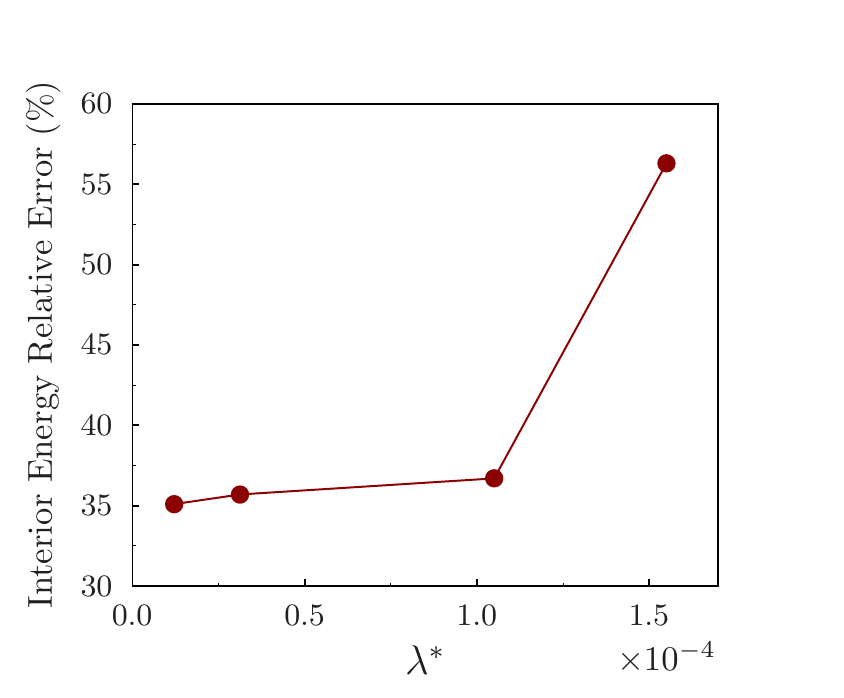}
   }

 \caption{Relation between the first discarded eigenvalue and the DG relative interior energy error;
               permeability from Fig.~\ref{fig:perm1} (left), permeability from Fig.~\ref{fig:perm2} (right)
               \ref{fig:perm1} and \ref{fig:perm2}
               }\label{fig:DGorder2}
\end{figure}

In the simulations using the DG formulation, we partition the original domain using a
coarse mesh of size $H=1/10$, and use a fine mesh composed of uniform triangular elements of mesh size $h=1/100$. The numerical results for permeability fields \ref{fig:perm1} and
\ref{fig:perm2} are represented in Tables \ref{table:dgperm1} and
\ref{table:dgperm2}, respectively. The first column shows the
dimension of the online space, the second column represents the
corresponding eigenvalue($\lambda^*$) of the first eigenfunction
discarded from the online space, and the next two columns illustrate the
interior energy relative error ($E_{\text{int}}$) and the boundary energy
relative error ($E_{\partial}$) between the fine scale solution and DG GMsFEM solution.
The errors between the offline and online solutions are offered in the final two columns.
We note that as the dimension of the online space increases (i.e., we keep more eigenfunctions in the space construction), the relative errors decrease accordingly. For example, the DG solution corresponding to Fig.~\ref{fig:perm1} yields interior relative energy errors that decrease from $55.08 - 34.86\%$, and boundary relative energy errors that decrease from $8.94 - 6.40 \%$. We note that in this case the Picard iteration converges in $4$ or $5$ steps for all simulations. In Fig.~\ref{fig:dgsolnperm2} we also plot the fine and coarse DG solutions that correspond to the field in Fig.~\ref{fig:perm2}. We note that the fine solution and the coarse solution corresponding to the smallest online space show some slight differences. However, the discrepancies noticeably diminish when the coarse DG solution is computed within the largest online space. As in the CG case, we also illustrate the relation between the DG interior errors and $\lambda^*$ in Fig.~\ref{fig:DGorder2}. From the plots in Fig.~\ref{fig:DGorder2}, we see that the relative errors decrease as $\lambda^*$ decreases, again following the expected error behavior.

\begin{table}
\tabcolsep 0pt %
\vspace*{1.5pt}
\begin{center}
\small \addtolength{\tabcolsep}{10.5pt}
\def\temptablewidth{1\textwidth}
\begin{tabular*}{\temptablewidth}{@{\extracolsep{\fill}}|c|c|c|c|c|c|}\hline \multirow{2}{*}{$\dim(V_{\text{on}}^{\text{DG}})$}
&\multirow{2}{*}{$\lambda^{*}$}&\multicolumn{2}{c|}{GMsFEM Relative Error (\%)}&\multicolumn{2}{c|}{Online-Offline Relative Error (\%)}\\
\cline{3-6} {}&{}&$\hspace{.6cm}E_{\text{int}}\hspace{.6cm}$ &$E_{\partial}$&

$\hspace{.8cm}E_{\text{int}}\hspace{.8cm}$  & $E_{\partial}$\\\hline\hline
       $271$& $1.53\times 10^{-4}$ &$55.08$&$8.94$&$44.38$&$8.43$ \\
\hline
       $331$&$1.24\times 10^{-4}$&$36.59$ &$6.63$ &$10.05$& $3.08 $ \\
\hline
       $466$&$3.03\times 10^{-5}$ &$35.57$&$6.56$ &$7.00$&$1.67$ \\
\hline
       $624$&$1.72\times 10^{-5}$ &$34.90$&$6.48$ &$2.12$&$0.40$ \\
\hline
       $716$&---& $34.86$&$6.40$&0.00&0.00\\
\hline\hline
       \end{tabular*}
       \end{center}
\caption{DG relative errors corresponding to the permeability field in
              Fig.~\ref{fig:perm1}; snapshot space uses $l_{\text{init}}^K$ eigenfunctions}
\label{table:dgperm1}
       \end{table}
%
%
\begin{table}
\tabcolsep 0pt %
\vspace*{1.5pt}
\begin{center}
\small \addtolength{\tabcolsep}{10.8pt}
\def\temptablewidth{1\textwidth}
\begin{tabular*}{\temptablewidth}{@{\extracolsep{\fill}}|c|c|c|c|c|c|}\hline \multirow{2}{*}{$\dim(V_{\text{on}}^{\text{DG}})$}
&\multirow{2}{*}{$\lambda^{*}$}&\multicolumn{2}{c|}{GMsFEM Relative Error (\%)}&\multicolumn{2}{c|}{Online-Offline Relative Error (\%)}\\
\cline{3-6} {}&{}&$\hspace{.6cm}E_{\text{int}}\hspace{.6cm}$ &$E_{\partial}$&

$\hspace{.7cm}E_{\text{int}}\hspace{.7cm}$  & $E_{\partial}$\\\hline\hline
       $270$& $1.56\times 10^{-4}$ &$56.29$&$10.30$&$46.37$&$9.75$ \\
\hline
       $331$&$1.05\times 10^{-4}$&$36.72$ &$6.71$ &$9.54$& $3.32$ \\
\hline
       $444$&$3.12\times 10^{-5}$ &$35.67$&$6.56$ &$6.48$&$1.67$ \\
\hline
       $582$&$1.21\times 10^{-5}$ &$35.06$&$6.48$ &$2.14$&$0.41$ \\
\hline
       $663$&---& $35.03$&$6.48$&0.00&0.00\\
\hline\hline
       \end{tabular*}
     \end{center}
\caption{DG relative errors corresponding to the permeability field in
              Fig.~\ref{fig:perm2}; snapshot space uses $l_{\text{init}}^K$ eigenfunctions}
\label{table:dgperm2}
\end{table}
%

\begin{remark}
When solving the nonlinear equation using the discontinuous Galerkin approach, we use different penalty parameters for fine-grid problem and coarse-grid problem (refer back to Subsect.~\ref{dgcoupling}). However, we observe that for different coarse penalty parameters that yield a convergent solution, the number of iterations and the relative errors (both interior and boundary) stay the same.
\end{remark}

\begin{remark}
Recall that we use the Galerkin projection of the previous coarse
solution onto the current online space as the approximation of the
previous coarse solution to obtain the terminal condition. If the coarse
penalty parameter is changed, we should use the current
coarse penalty parameter to construct the Galerkin projection.
\end{remark}

We observe from Tables \ref{table:cgperm1}-\ref{table:dgperm2} that
the offline spaces for DG formulation are much smaller than those
obtained through CG formulation. As a result, in Table \ref{table:extdgperm2}
we use more eigenfunctions (more specifically, we set $l_{\text{max}}^K = l_{\text{init}}^K+3$) in the snapshot space construction to yield a larger
offline space. For these examples, we use the permeability field from Fig.~\ref{fig:perm2}. Due to the increase of the offline (and corresponding online) space dimensions, we see more accurate results than those offered in Table \ref{table:dgperm2}.

\begin{table}
\tabcolsep 0pt %
\vspace*{1.5pt}
\begin{center}
\small \addtolength{\tabcolsep}{10.8pt}
\def\temptablewidth{1\textwidth}
\begin{tabular*}{\temptablewidth}{@{\extracolsep{\fill}}|c|c|c|c|c|c|}\hline \multirow{2}{*}{$\dim(V_{\text{on}}^{\text{DG}})$}
&\multirow{2}{*}{$\lambda^{*}$}&\multicolumn{2}{c|}{GMsFEM Relative Error (\%)}&\multicolumn{2}{c|}{Online-Offline Relative Error (\%)}\\
\cline{3-6} {}&{}&$\hspace{.6cm}E_{\text{int}}\hspace{.6cm}$ &$E_{\partial}$&

$\hspace{.8cm}E_{\text{int}}\hspace{.8cm}$  & $E_{\partial}$\\\hline\hline
       $381$& $1.47\times 10^{-4}$ &$37.34$&$7.42$&$22.80$&$6.00$ \\
\hline
       $440$&$1.54\times 10^{-4}$&$35.92$ &$6.16$ &$20.07$& $4.36 $ \\
\hline
       $707$&$9.54\times 10^{-5}$ &$32.80$&$5.29$ &$13.64$&$2.90$ \\
\hline
       $958$&$2.71\times 10^{-5}$ &$29.44$&$5.48$ &$4.80$&$0.98$ \\
\hline
       $1352$&---& $28.98$&$5.39$&0.00&0.00\\
\hline\hline
       \end{tabular*}

      \end{center}
\caption{DG relative errors corresponding to the permeability field in
              Fig.\ref{fig:perm2}; snapshot space uses $l_{\text{max}}^K = l_{\text{init}}^K+3$ eigenfunctions}
        \label{table:extdgperm2}
\end{table}

\subsection{Parameter-dependent permeability field}
\label{withparameter}
For the next set of numerical results, we consider solving the nonlinear elliptic problem
in Eq.~\eqref{eqn:problem}
with a coefficient of the form $\kappa(x,u,\mu^{\text{p}})=\exp\left[\left(\mu^{\text{p}}
\kappa_{1}(x)+(1-\mu^{\text{p}})\kappa_{2}(x)\right)u(x)\right]$. For $\kappa_1(x)$ and $\kappa_2(x)$ we use the fields shown in Fig. \ref{fig:permi} and \ref{fig:permii}, respectively. As for the parameter-dependent simulation, we are careful to distinguish the difference between the auxiliary parameter $\mu = \overline{u}^n$ which is used to denote a previous solution iterate, and a ``physical" parameter $\mu^{\text{p}}$ that is used in the construction of a new permeability field. We take the range of $\mu^{\text{p}}$ to be $[0, 1]$, and use three equally spaced values in order to construct the snapshot space in this case. We use the same $[u_{\text{min}}, u_{\text{max}}]$ interval from the previous results, yet use four equally spaced values in this case. In particular, we use the pairs $(u_j,\mu_l^{\text{p}})$, where $1\leq j \leq 4$, and $1\leq l \leq 3$ as the fixed parameter values for the snapshot space construction. At the online stage we use the initial guess $u^0=0$ and a fixed value of $\mu^{\text{p}} = 0.2$ while solving the respective eigenvalue problem required for the continuous or discontinuous Galerkin online space construction.

In Table~\ref{table:paracgperm} we offer results corresponding to the CG formulation, and in Tables~\ref{table:paradgperm} and \ref{table:extparadgperm} we offer results corresponding to the DG formulation. In all cases we encounter very similar error behavior compared to the examples offered earlier in the section. In particular, an increase of the dimension of the online space yields predictably smaller errors, and smaller values of $\lambda^*$ correspond to the error decrease. And while it suffices to refer back to related discussions earlier in the section, we emphasize that this distinct set of results serves to further illustrate the robustness of the proposed method. In particular, we show that the solution procedure allows for a suitable treatment of nonlinear problems that involve auxiliary and physical parameters.

\begin{table}
\tabcolsep 0pt %
\vspace*{1.5pt}
\begin{center}
\small \addtolength{\tabcolsep}{8.9pt}
\def\temptablewidth{1\textwidth}
\begin{tabular*}{\temptablewidth}{|c|c|c|c|c|c|}
\hline \multirow{2}{*}{$\dim(V_{\text{on}}^{\text{CG}})$}
&\multirow{2}{*}{$\lambda^{*}$}&
\multicolumn{2}{c|}{\;\;\;GMsFEM Relative Error (\%)\;\;\;}&\multicolumn{2}{c|}{\hspace{0cm}Online-Offline Relative Error (\%)\hspace{-.5cm}}\\
\cline{3-6} {}&{}&\hspace{.6cm}$L^{2}_{\kappa}(D)\hspace{.6cm}$ &
$H^{1}_{\kappa}(D)$&$\hspace{.6cm}L^{2}_{\kappa}(D)\hspace{.6cm}$&
$H^{1}_{\kappa}(D)$
\\
\hline\hline       
       $309$&$0.0027$ &$1.30$&$14.89$ &$1.10$&$15.32$ \\

\hline
      $492$&$0.0010$ & $0.59$&$10.82$& $0.39$&$9.76$ \\
\hline
       $580$&$6.76\times 10^{-4}$& $0.45$&$9.55$ & $0.24$&$7.92$\\
\hline
      $728$&$3.33\times 10^{-4}$&$0.34$&$7.87$ &$0.12$&$5.23$ \\
\hline
       $991$& ---&$0.28$&$6.74$&0.00&0.00\\
\hline\hline
\end{tabular*}
      \end{center}
\caption{CG relative errors corresponding to the parameter-dependent field constructed from
              Fig.~\ref{fig:permi} and \ref{fig:permii} }
\label{table:paracgperm}
       \end{table}
%
\begin{table}
\tabcolsep 0pt %
\vspace*{1.5pt}
\begin{center}
\small \addtolength{\tabcolsep}{10.8pt}
\def\temptablewidth{1\textwidth}
\begin{tabular*}{\temptablewidth}{@{\extracolsep{\fill}}|c|c|c|c|c|c|}\hline \multirow{2}{*}{$\dim(V_{\text{on}}^{\text{DG}})$}
&\multirow{2}{*}{$\lambda^{*}$}&\multicolumn{2}{c|}{GMsFEM Relative Error (\%)}&\multicolumn{2}{c|}{Online-Offline Relative Error (\%)}\\
\cline{3-6} {}&{}&$\hspace{.6cm}E_{\text{int}}\hspace{.6cm}$ &$E_{\partial}$&

$\hspace{.8cm}E_{\text{int}}\hspace{.8cm}$  & $E_{\partial}$\\\hline\hline
       $300$& $1.02\times 10^{-4}$ &$37.56$&$7.94$&$10.15$&$3.16$ \\
\hline
       $313$&$6.25\times 10^{-5}$&$37.55$ &$7.81$ &$10.00$& $2.85$ \\
\hline
       $403$&$2.58\times 10^{-5}$ &$36.81$&$7.35$ &$5.83$&$1.38$ \\
\hline
       $497$&$1.22\times 10^{-5}$ &$36.37$&$7.21$ &$0.84$&$0.10$ \\
\hline
       $517$&---& $36.36$&$7.21$&0.00&0.00\\
\hline\hline
       \end{tabular*}

       \end{center}
\caption{DG relative errors corresponding to the parameter-dependent field constructed from
              Fig.~\ref{fig:permi} and \ref{fig:permii}; snapshot space uses $l_{\text{init}}^K$ eigenfunctions}
\label{table:paradgperm}
\end{table}
%
\begin{table}
\tabcolsep 0pt %
\vspace*{1.5pt}
\begin{center}
\small \addtolength{\tabcolsep}{10.8pt}
\def\temptablewidth{1\textwidth}
\begin{tabular*}{\temptablewidth}{@{\extracolsep{\fill}}|c|c|c|c|c|c|}\hline \multirow{2}{*}{$\dim(V_{\text{on}}^{\text{DG}})$}
&\multirow{2}{*}{$\lambda^{*}$}&\multicolumn{2}{c|}{GMsFEM Relative Error (\%)}&\multicolumn{2}{c|}{\hspace{-.2cm}Online-Offline Relative Error (\%)\hspace{-.2cm}}\\
\cline{3-6} {}&{}&$\hspace{.6cm}E_{\text{int}}\hspace{.6cm}$ &$E_{\partial}$&

$\hspace{.8cm}E_{\text{int}}\hspace{.8cm}$  & $E_{\partial}$\\\hline\hline
       $300$& $2.13\times 10^{-4}$ &$37.59$&$7.94$&$22.54$&$6.40$ \\
\hline
       $440$&$1.54\times 10^{-5}$&$35.78$ &$5.92$ &$18.89$& $3.74$ \\
\hline
       $668$&$7.69\times 10^{-5}$ &$32.54$&$5.39$ &$11.62$&$2.58$ \\
\hline
       $902$&$1.51\times 10^{-5}$ &$30.23$&$5.29$ &$3.87$&$1.06$ \\
\hline
       $1093$&---& $29.88$&$5.29$&0.00&0.00\\
\hline\hline
       \end{tabular*}

       \end{center}
\caption{DG relative errors corresponding to the parameter-dependent field constructed from
              Fig.~\ref{fig:permi} and \ref{fig:permii}; snapshot space uses $l_{\text{max}}^K = l_{\text{init}}^K+3$ eigenfunctions}
\label{table:extparadgperm}
\end{table}
%

\section{Concluding Remarks}
\label{conclusion}
In this paper we use the Generalized Multiscale Finite Element (GMsFEM) framework in order to solve nonlinear elliptic equations with high-contrast coefficients. In order to solve this type of problem we linearize the equation such that upscaled quantities of previous solution iterates may be regarded as auxiliary coefficient parameters in the problem formulation. As a result, we are able to construct a respective set of coarse basis functions using an offline-online procedure in which the precomputed offline space allows for the efficient computation of a smaller-dimensional online space for any parameter value at each iteration. In this paper, the coarse space construction involves solving a set of localized eigenvalue problems that are tailored to either continuous Galerkin (CG) or discontinuous Galerkin (DG) global coupling mechanisms. In particular, the respective coarse spaces are formed by keeping a set of eigenfunctions that correspond to the localized eigenvalue behavior. Using either formulation, we show that the process of systematically enriching the coarse solution spaces yields a predictable error decline between the fine and coarse-grid solutions. As a result, the proposed methodology is shown to be an effective and flexible approach for solving the nonlinear, high-contrast elliptic equation that we consider in this paper.

\section*{Acknowledgements}

Y. Efendiev's work is
partially supported by the
DOE and NSF (DMS 0934837 and DMS 0811180).
J.Galvis would like to acknowledge partial support from DOE.

 This publication is based in part on work supported by Award
No. KUS-C1-016-04, made by King Abdullah University of Science
and Technology (KAUST).

\section*{References}

\bibliographystyle{plain}

\def\cprime{$'$}

\end{document}